\documentclass[10pt]{article}
\usepackage[francais,english]{babel}
\usepackage{amsmath}
\usepackage{amsfonts}
\usepackage{amssymb}
\usepackage{amsthm}
\usepackage{graphics}
\usepackage{graphicx}
\usepackage{amscd}
\usepackage{fullpage}
\usepackage{epsfig}
\usepackage{color}
\usepackage[matrix,arrow,ps]{xy}

\newcommand{\Z}{\mathbb{Z}}
\newcommand{\R}{\mathbb{R}}
\newcommand{\N}{\mathbb{N}}

\newcommand{\C}{\mathbb{C}}

\newcommand{\E}{\mathbb{E}}

\newcommand{\mc}{\mathcal}
\newcommand{\mb}{\mathbb}
\newcommand{\mf}{\mathfrak}
\newcommand{\eps}{\varepsilon}
\newcommand{\ind}{{\bf 1}}

\newcommand{\rK}{{\sf K}}

\newcommand{\dg}{\dagger}
\newcommand{\dmd}{\diamondsuit}

\newcommand{\Mon}{{\sf Mon}}

\renewcommand{\H}{\mathbb{H}}

\newcommand{\Lap}{\Delta\!}

\DeclareMathOperator{\Pf}{Pf}

\title{Exact bosonization of the Ising model}
\author{Julien Dub\'edat\footnote{Partially supported by NSF grant DMS-1005749 and the Alfred P. Sloan Foundation.}}
\newtheorem{thm}{Theorem}

\newtheorem{Thm}[thm]{Theorem}

\newtheorem{Lem}[thm]{Lemma}

\begin{document}
\maketitle
\begin{abstract}
We present exact combinatorial versions of bosonization identities, which equate the product of two Ising correlators with a free field (bosonic) correlator. The role of the discrete free field is played by the height function of an associated bipartite dimer model. Some applications to the asymptotic analysis of Ising correlators are discussed.
\end{abstract}

\section{Introduction}

The Ising model is a basic model of ferromagnetism in statistical mechanics. It has been extensively studied in two dimensions, at and around its critical temperature (\cite{MW_Ising}), where several powerful techniques are available: transfer matrices, pfaffian (dimer) representation, integrable systems and more recently discrete complex analysis and Schramm-Loewner Evolutions (\cite{Smi_ICM}). 

In the Ising model, the basic (order) observables are the spin variables; Kramers-Wannier duality (\cite{ KW_ising1,KW_ising2}) maps these to disorder variables (\cite{KC_disorder}). In general, an order variable is a local random variable, while a disorder variable represents a local modification of the state space (eg \cite{Dub_abelian}). A classical object of study are correlators involving some order and/or disorder variables, especially in appropriate asymptotic regimes (\cite{MW_Ising,Palmer_planar}).

At the field theoretic level, the notion of bosonization introduced in \cite{ZubItz_Ising} and much developed afterwards (see in particular \cite{DiFSalZub_isingtorus} and Chapter 12 in \cite{DiF}) expresses squares of Ising correlators as bosonic (free field) correlators. For the free field, the relevant order variables are electric insertions, and the  dual disorder variables are magnetic insertions (eg \cite{DiF,Gaw_CFT,Dub_abelian}).

It has long been known that a planar Ising configuration may be represented by a dimer configuration on a related decorated graph (\cite{Fisher_Ising,Kas_Ising}), which may be analyzed through determinantal or Pfaffian techniques (\cite{Kas_square}); these decorated graphs are not bipartite. Dimers on bipartite graphs are associated to a discrete height function, the asymptotic fluctuations of which have been extensively studied from the early 90's (see the survey \cite{Ken_IAS}). In the type of scaling regime relevant here, these asymptotic fluctuations are described by a free field in a rather precise fashion.

In the present article, we use (known) mappings and dualities between doubled Ising models, 8-vertex, 6-vertex and dimer models (along the lines of \cite{Baxter_exact}), tracking order and disorder variables along the way, in order to establish exact combinatorial versions of the field-theoretic bosonization identities. These identities involve a pair of independent Ising configurations and a single bipartite dimer configuration. In combination with recent progress on the fine limiting behavior of dimer height fields (\cite{Dub_tors}), this enables to obtain the asymptotics of critical Ising correlators in the plane. Critical correlators in finite domains require some additional care and are the object of the upcoming \cite{CHI,Dub_prep}.

The article is organized as follows. Bosonization identities are phrased in Sect. 2, along with a discussion of boundary conditions. Some consequences for the asymptotic analysis of critical Ising correlators are listed in Sect. 3. Relations with other approaches are described in Sect. 4.

\section{Bosonization rules}

\subsection{Mappings}

Consider a graph $\Gamma=(V,E)$ embedded on a torus $\Sigma=\C/\Lambda$, $\Lambda=\Z+\tau\Z$, $\Im\tau>0$ (the planar case will be discussed afterwards). Let $F$ denote the set of its faces and $\Gamma^\dg=(V^\dg,E^\dg)$ denote its dual graph, so that $V^\dg\simeq F$, $E^\dg\simeq E$ (if $e\in E$, $e^\dg\in E^\dg$ denotes its dual edge). To each edge $e\in E$ we associate a coupling constant $J_e\geq 0$. A configuration of the Ising model consists in an assignment of a spin $\sigma_v=\pm 1$ to each vertex $v$ of $\Gamma$. By planar duality, one can equivalently assign spins to faces. The weight of a configuration $(\sigma_v)_{v\in V}$ is:
$$w((\sigma_v)_{v\in V})=\exp(-2\beta \sum_{e=(vv')\in E} J_e\ind_{\sigma_v\neq\sigma_{v'}})$$
where the inverse temperature $\beta$ is a fixed positive constant. The weight is invariant under global spin flip $(\sigma_v)_v\leftrightarrow (-\sigma_v)_v$. In the {\em low temperature expansion}, one represents a configuration $(\sigma_v)$ by the even degree subgraph ({\em polygon}) $P_\sigma=(V^\dg,E^\dg_\sigma)$ of $\Gamma^\dg=(V^\dg,E^\dg)$, where $e=(vv')^\dg\in E_\sigma$ iff $\sigma_v\neq\sigma_{v'}$. Clearly
$$w((\sigma_v))=\prod_{e\in E^\dg_\sigma} w(e)$$
if we set 
$$w(e)=\exp(-2\beta J_{e^\dg})\in (0,1]$$
A polygon $P$ has even degree at each vertex. Given $P$, it is easy to see that the parity of the number of edges crossed by a closed cycle on $\Gamma$ depends only on the homology class of this cycle (in $H_1(\Sigma,\Z)$). Thus to any polygon one may associate signs $(\eps_A(P),\eps_B(P))\in\{\pm 1\}^2$ where $\eps_A=+1$ or $-1$ according to whether an $A$-cycle crosses an even or odd number of edges of $P$, and similarly for $\sigma_B$. Plainly, a polygon $P$ comes from a (periodic) spin configuration $\sigma$ iff $(\eps_A(P),\eps_B(P))=(1,1)$. (We could frame this discussion in terms of discrete $1$-forms with values in $\{\pm 1\}$).

At this point it is natural to introduce spin configurations with periodic ($p$) or antiperiodic ($a$) boundary conditions corresponding to $A$ and $B$ periods. One way is to consider spin configurations $\sigma$ on the lift of $\Gamma$ to $\C/(2\Lambda)$ such that $\sigma(x+1)=\eps_A\sigma(x)$, $\sigma(x+\tau)=\eps_B\sigma(x)$, where $(\eps_A,\eps_B)=(1,1),(-1,1),(1,-1),(-1,-1)$ corresponds to $(pp),(ap),(pa),(aa)$ boundary conditions. Associated to such a spin configuration $\sigma$, we have a polygon $P_\sigma$, where now $(\eps_A,\eps_B)(P_\sigma)$ may take any value in $\{\pm 1\}^2$.

Thus consider the space of polygons on $\Gamma$ with weights
$$w(P)=\prod_{e\in E_P}w(e)$$
The state space is partitioned in four blocks corresponding to $(pp)$,\dots,$(aa)$ boundary condition. Each polygon lifts to two spin configurations $\pm\sigma$, in a measure preserving way (up to normalization).

Let us also introduce the {\em medial graph} $M=(V_M,E_M)$ of $\Gamma$ (or derived graph, in the restricted context of planar graphs): the vertices $V_M$ of $M$ are set at the midpoints of edges of $\Gamma$; two vertices of $M$ are adjacent if they correspond to consecutive (in cyclic order) edges around a vertex of $\Gamma$. Note that $M$ is $4$-regular (all vertices have degree 4) and may be identified with the medial graph of $\Gamma^\dg$. Every edge $e$ of $M$ corresponds to two edges of $\Gamma$ that have a common endpoint $v\in V$ and lie on the boundary of a common face $f\in F$. This establishes a natural correspondence between edges of $M$ and pairs $(v,f)\in V\times F$ where $v$ is a vertex on the boundary of $F$.

\begin{figure}[htb]
\begin{center}
\leavevmode
\includegraphics[width=0.8\textwidth]{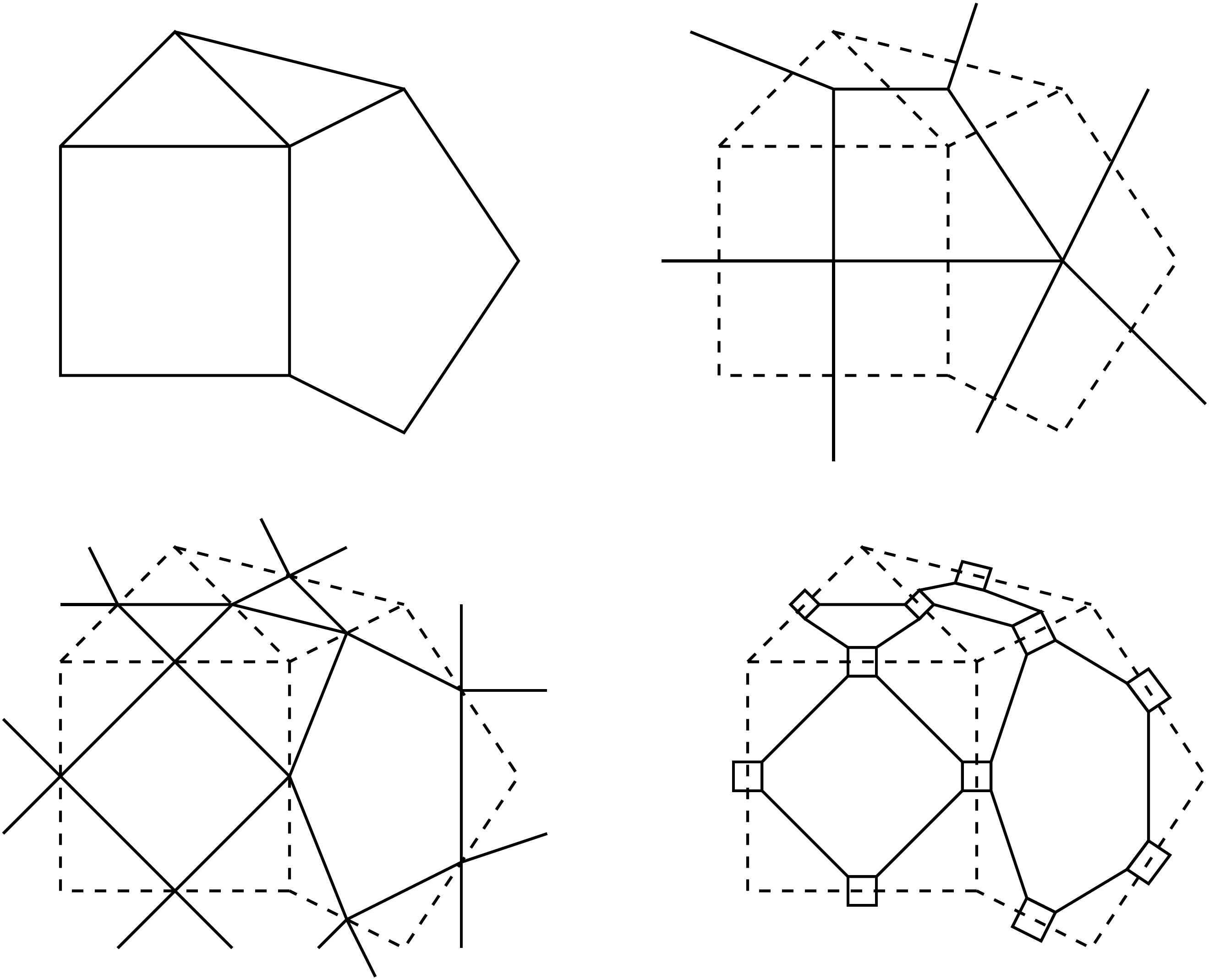}
\end{center}
\caption{(portion of a) graph $\Gamma$; dual graph $\Gamma^\dg$; medial graph $M$; bipartite graph $C$}
\label{fig:medial}
\end{figure}

Given a $4$-regular planar graph (such as $M$), we define an {\em $8$-vertex configuration} (or $8V$ configuration for short) as an orientation of the graph such that the number of incoming (or outgoing) edges at each vertex is even. Around each vertex, there are 8 possible configurations (orientations of the four adjacent edges).

{\bf From dual Ising models to the 8V model.}

\begin{figure}[htb]
\begin{center}
\leavevmode
\includegraphics[width=0.8\textwidth]{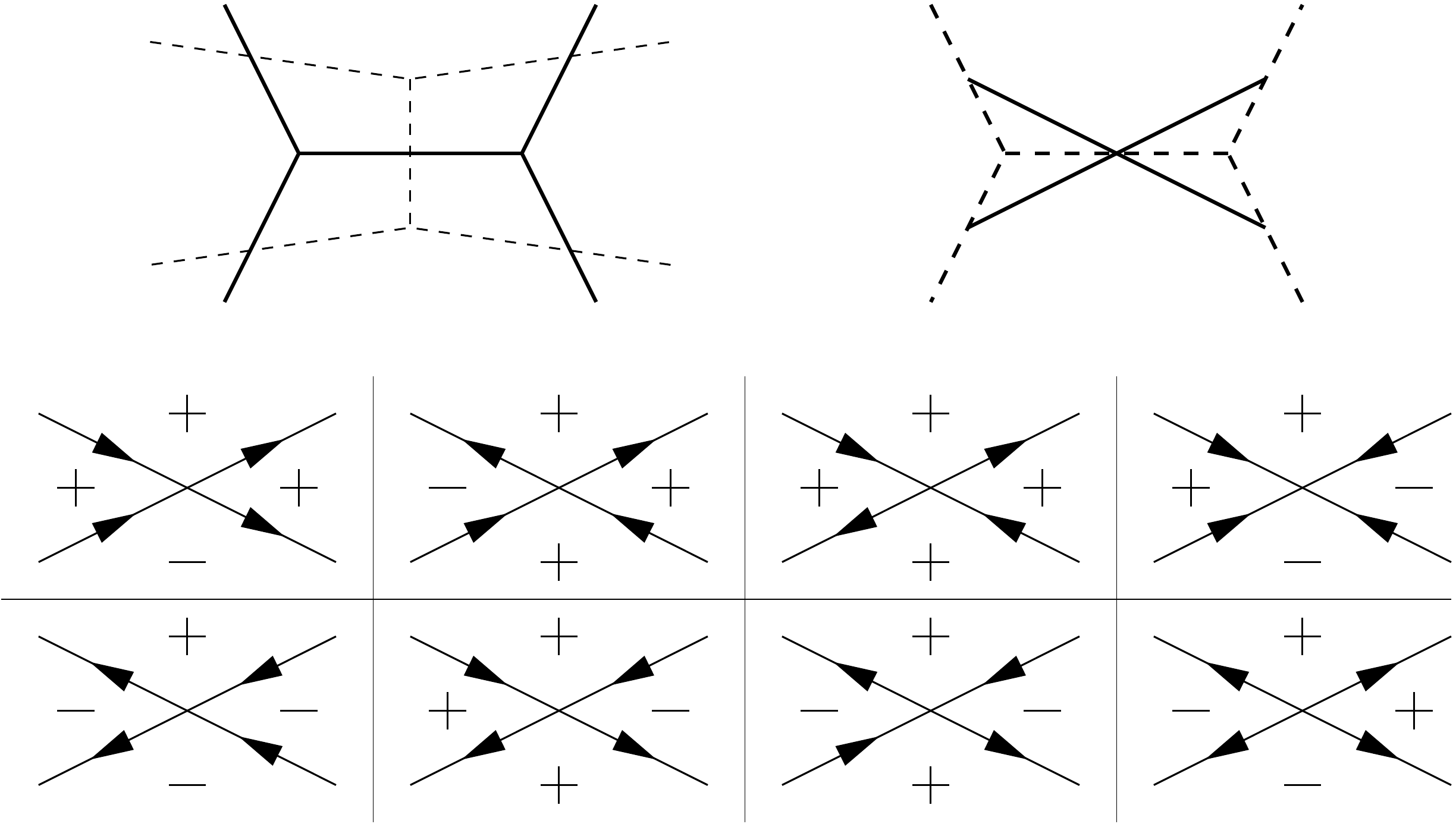}
\end{center}
\caption[Graph $\Gamma$ (solid) and $\Gamma^\dg$ (dashed); medial graph; spin variables and associated 8V configuration type:]{Graph $\Gamma$ (solid) and $\Gamma^\dg$ (dashed); medial graph; spin variables and associated 8V configuration type:\begin{tabular}{c|c|c|c}1&3&5&7\\ 2&4&6&8\end{tabular}}
%
%
\label{fig:8V}
\end{figure}

Let us now consider two spin configurations $\sigma$ and $\sigma^\dg$, with $\sigma$ defined on $V$ and $\sigma^\dg$ defined on $F\simeq V^\dg$ with the same boundary condition in $\{(pp),\dots,(aa)\}$. 
To the pair $(\sigma,\sigma^\dg)$ we associate an $8V$ configuration on $M$ as follows: Let $e$ be an edge of $M$ corresponding to the pair $(v,f)\in V\times F$. If $\sigma(v)=\sigma^\dg(f)$ (resp. $\sigma(v)=-\sigma^\dg(f)$), we orient $e$ such that $v$ is on its right handside (resp. left handside). In other words we have a reference orientation of $M$ where vertices of $\Gamma$ are on the righthand side of oriented edges of $M$. Other orientations are associated with an edge spin configuration: $\nu(e)=\sigma(v)\sigma^\dg(f)$, where $e\in E_M$ corresponds to $(v,f)\in V\times F$ (see Figure \ref{fig:8V}). The orientation (8V configuration) associated to $\nu$ agrees with the reference orientation on $e$ iff $\nu(e)=1$.

Let us now assume that $(\sigma(v))_{v\in V}$ and $(\sigma^\dg(f))_{f\in F}$ are sampled independently from Ising distributions (with couplings $(J_e)_{e\in E}$, $(J_{e^\dg})_{e^\dg\in E^\dg}$ and associated edge weights $(w(e))_e$, $(w(e^\dg))_{e^\dg}$) and same boundary condition in $\{(pp),\dots,(aa)\}$. The weight of the configuration $(\sigma,\sigma^\dg)$ (up to normalization) may be expressed as a product of factors $w(e)w(e^\dg)$, where $e,e^\dg$ is a pair of dual edges of $\Gamma,\Gamma^\dg$ that corresponds to a vertex of $M$. The 8V weights at this vertex are as follows ($\omega_i$ is the weight of the $i$-th vertex configuration, see Figure \ref{fig:8V}):
\begin{eqnarray*}
\omega_1=\omega_2=w(e)&\omega_5=\omega_6=1\\
\omega_3=\omega_4=w(e^\dg)&\omega_7=\omega_8=w(e)w(e^\dg)
\end{eqnarray*}
Remark that $(\sigma,\sigma^\dg)$ and $(-\sigma,-\sigma^\dg)$ project to the same 8V configuration, while $(\sigma,-\sigma^\dg)$, $(-\sigma,\sigma^\dg)$ projects to the configuration with all arrows reverted (it has the same weight). When 8V weights are invariant under reversal of all arrows, we denote $\omega_{12}=\omega_1=\omega_2$ etc.

{\bf From the 8V model to the 6V model.}

These weights may be written as a function of the medial edge spins $\nu(e_m)=\sigma(v)\sigma^\dg(f)$. If we number \begin{tabular}{cc} 1&4\\ 2&3\end{tabular} the medial edges in Figure \ref{fig:8V}, we may write the weight as
$$\omega(\nu_1,\dots,\nu_4)=a(1+\nu_1\nu_2\nu_3\nu_4)+b(\nu_1\nu_2+\nu_3\nu_4)+c(\nu_1\nu_4+\nu_2\nu_3)+d(\nu_1\nu_3+\nu_2\nu_4)$$
because of the invariances $\omega(\nu_1,\dots,\nu_4)=\omega(-\nu_1,\dots,-\nu_4)$ and $\omega(\nu_1,\nu_2,\nu_3,\nu_4)=\omega(\nu_3,\nu_4,\nu_1,\nu_2)$. Note that $\omega(\nu_1,\dots,\nu_4)=0$ if $\nu_1\dots\nu_4=-1$, which enforces the 8V condition. Solving
$$\left\{\begin{array}{ll}
2(a-b+c-d)&=\omega(1,-1,-1,1)=\omega_{12}\\
2(a+b-c-d)&=\omega(1,1,-1,-1)=\omega_{34}\\
2(a+b+c+d)&=\omega(1,1,1,1)=\omega_{56}\\
2(a-b-c+d)&=\omega(1,-1,1,-1)=\omega_{78}
\end{array}\right.$$
we get
$$8(a,b,c,d)=(\omega_{12}+\omega_{34}+\omega_{56}+\omega_{78},-\omega_{12}+\omega_{34}+\omega_{56}-\omega_{78},
\omega_{12}-\omega_{34}+\omega_{56}-\omega_{78},
-\omega_{12}-\omega_{34}+\omega_{56}+\omega_{78}
)$$

The abelian duality for the 8-vertex model is obtained as follows. Let us start from the 8V partition function
$${\mc Z}=\sum_{(\nu_e)\in\{\pm 1\}^{E_M}}\prod_{e\in E}\omega_e(\nu)$$
where $\omega_e(\nu)$ is the $8V$ weight at $e$, which is written as a sum of characters of $\{\pm 1\}^4$. Expanding $\prod_e\omega_e(\nu)$ in monomials in the $\nu_e$ variables, we notice that only monomials with even degree (0 or 2) in each variable contribute to the partition function. To each contributing monomial, we associate a medial spin configuration $(\hat \nu_e)_{e\in E_M}$ as follows: $\hat\nu_e=1$ if the partial degree of $\nu_e$ is 0 and $\hat\nu_e=-1$ if the partial degree of $\nu_e$ is 2. This yields
$${\mc Z}=\hat{\mc Z}$$
where $\hat{\mc Z}$ is the 8V partition function with weights
\begin{align*}
\hat\omega_{12}&=c/2=\frac{\omega_{12}-\omega_{34}+\omega_{56}-\omega_{78}}4\\
\hat\omega_{34}&=b/2=\frac{-\omega_{12}+\omega_{34}+\omega_{56}-\omega_{78}}4\\
\hat\omega_{56}&=a/2=\frac{\omega_{12}+\omega_{34}+\omega_{56}+\omega_{78}}4\\
\hat\omega_{78}&=d/2=\frac{-\omega_{12}-\omega_{34}+\omega_{56}+\omega_{78}}4
\end{align*}
Note that the weight mapping $(\omega_{12},\dots)\mapsto (\hat\omega_{12},\dots)$ is involutive. 
The spin variables $\nu_e$, $e\in E_M$, are natural order variables. It is easy to see that under duality, they are exchanged with disorder variables defined as follows. A disorder variable $\xi_e$ at the medial edge $e\in E_M$ is a defect splitting it in two half-edges with opposite orientations. It is thus a modification of the state space, which we still denote as a random variable by (a standard) abuse of terminology. The duality identity then reads:
$$\langle \nu(e_1)\dots\nu(e_m)\xi(e_{m+1})\dots\xi(e_{m+n})\rangle_\omega= \langle \xi(e_1)\dots\xi(e_m)\nu(e_{m+1})\dots\nu(e_{m+n})\rangle_{\hat\omega}$$
where the LHS is $\sum \nu(e_1)\dots\nu(e_m)\prod_{e\in E}\omega_e(\nu)$, where the sum bears on configurations with defects at $e_{m+1},\dots,e_{m+n}$. The RHS is defined symmetrically. 

In presence of a disorder $\xi_e$, there are two opposite spin variables $\nu(e^-)$, $\nu(e^+)$ corresponding to the two half-edges $e^+,e^-$ of $e$. Then duality maps the pair $\nu(e^+)\xi(e)$ to $\xi(e)\nu(e^-)$, and vice-versa.

Assume that the two sets of Ising weights $(w(e))_{e\in E}$, $(w(e^\dg))_{e^\dg\in E^\dg}$ satisfy the Kramers-Wannier duality relation:
\begin{equation}\label{eq:KW}
w(e)+w(e^\dg)+w(e)w(e^\dg)=1
\end{equation}
for each pair $(e,e^\dg)$ of dual edges. Then the associated 8V weights are
\begin{equation}\label{eq:8Vweights}
(\omega_{12},\omega_{34},\omega_{56},\omega_{78})=(w,w',1,ww')
\end{equation}
with $w=w(e)$, $w'=w(e^\dg)$ for short. In any 8V configuration on a toroidal graph, the number of sinks (type 7) equals the number of sources (type 8). Consequently, we get the same configuration weights if we change the local weights to:
\begin{equation}\label{eq:8Vweights2}
(\omega_{12},\omega_{34},\omega_{56},\omega_{78})=(w,w',1,-ww')
\end{equation}
Applying duality, we get the same partition function with the dual weights:
\begin{align*}
(\hat\omega_{12},\hat\omega_{34},\hat\omega_{56},\hat\omega_{78})&=\frac 14(w-w'+1+ww',-w+w'+1+ww',w+w'+1-ww',-w-w'+1-ww')\\
&=\frac 12(1-w',1-w,1-ww',0)
\end{align*}
This set of weights defines a 6-vertex model (as sources and sinks get zero weight), with weights:
\begin{equation}\label{eq:6Vweights}
(\hat\omega_{12},\hat\omega_{34},\hat\omega_{56})
=\frac 12(1-w',1-w,1-ww')
\end{equation}
Moreover,
$$\frac{\hat\omega_{12}}{\hat\omega_{56}}=\frac{2w}{1+w^2},\ \ \ 
\frac{\hat\omega_{34}}{\hat\omega_{56}}=\frac{1-w^2}{1+w^2}$$
so that $\hat\omega_{56}^2=\hat\omega_{12}^2+\hat\omega_{34}^2$, ie these are the weights of a 6V model at the ``free fermion" point.

{\bf From 6V to bipartite dimers.}

Fan and Wu showed (\cite{FanWu}) that the 8V model on the free fermion line, ie satisfying the vertex weight relation:
$$\omega_1\omega_2+\omega_3\omega_4=\omega_5\omega_6+\omega_7\omega_8$$
can be mapped exactly to a dimer model on a decorated graph, in a way rather similar to the Temperley-Fisher mapping of the Ising model to a dimer model. In the case of the 6V model at the free fermion point, one can find a bipartite dimer representation, which we now describe.

We start from a 4-regular planar graph $M$ (plainly, this works for toroidal graphs), with checkerboard coloring of faces. A decorated graph $C$ is obtained as follows. Each vertex $v$ of $M$ is replaced with a {\em city} (as in the urban renewal metaphor, see \cite{KPW}), a quadrangle with a vertex on each edge abutting $v$. Each city has four internal edges and two adjacent cities are connected by a road edge. Formally, vertices of $C$ correspond to pairs $(v,e)$, $v\in M$, $e\in E_M$ abutting $v$; $((v,e),(v',e'))$ is an edge of $C$ iff either $e=e'=(vv')$ (road) or $v=v'$ and $e,e'$ are two consecutive (in cyclic order) edges around $v$ (city street).

Plainly, $C$ is bipartite: the city faces are quadrangles, and the degree of other faces is twice the degree of the corresponding face in $M$ (an additional argument is needed to ensure consistency when wrapping around the torus; this is obtained from the checkerboard coloring of the faces of $M$). A {\em dimer configuration} or {\em perfect matching} on $C$ is a subset ${\mf m}$ of edges of $C$ such that each vertex of $C$ is the endpoint of exactly one edge in $C$. Given a set $(w_e)_{e\in E_C}$ of positive weights associated to edges of $C$, one defines a dimer configuration weight by
$$w({\mf m})=\prod_{e\in{\mf m}}w_e$$
For background on the dimer model, see eg \cite{Ken_IAS}.

\begin{figure}[htb]
\begin{center}
\leavevmode
\includegraphics[width=0.8\textwidth]{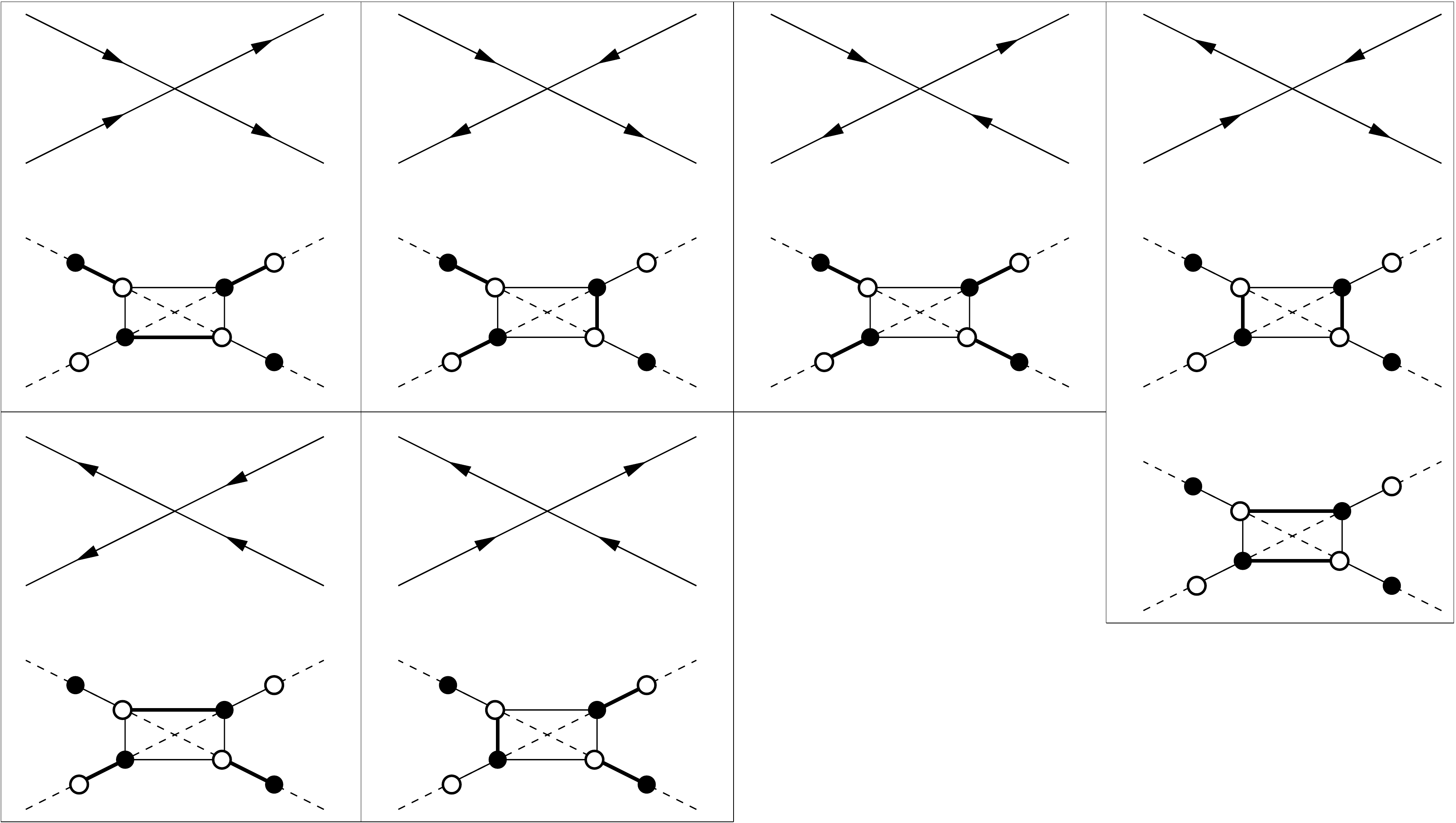}
\end{center}
\caption{local 6V configurations on a 4-regular graph (type 1,2;3,4;5;6) and associated local dimer configurations on a decorated graph (bold edge: dimer)
}
%
%
\label{fig:6Vcity}
\end{figure}

We can choose the checkerboard coloring of faces of $M$ and the two-coloring of vertices of $C$ in such a way that if $(bw)$ is a road edge of $C$ oriented from black to white, the face of $M$ on the RHS (resp. LHS) is black (resp. white). A 6V configuration consists on a orientation of edges of $M$ (with the 6V rule enforced at each vertex). Given such a configuration, each road edge of $C$ inherits an orientation. We decide that $(bw)\in E_C$ belongs to the corresponding dimer configuration iff $(bw)$ is oriented from black to white. This completely determines the perfect matching, up to a local ambiguity for each type 6 city, where two opposite city streets are matched. See Figure \ref{fig:6Vcity} for a graphical representation of the correspondence.

\begin{figure}[htb]
\begin{center}
\leavevmode
\scalebox{1}[1.73]{\includegraphics[width=.8\textwidth]{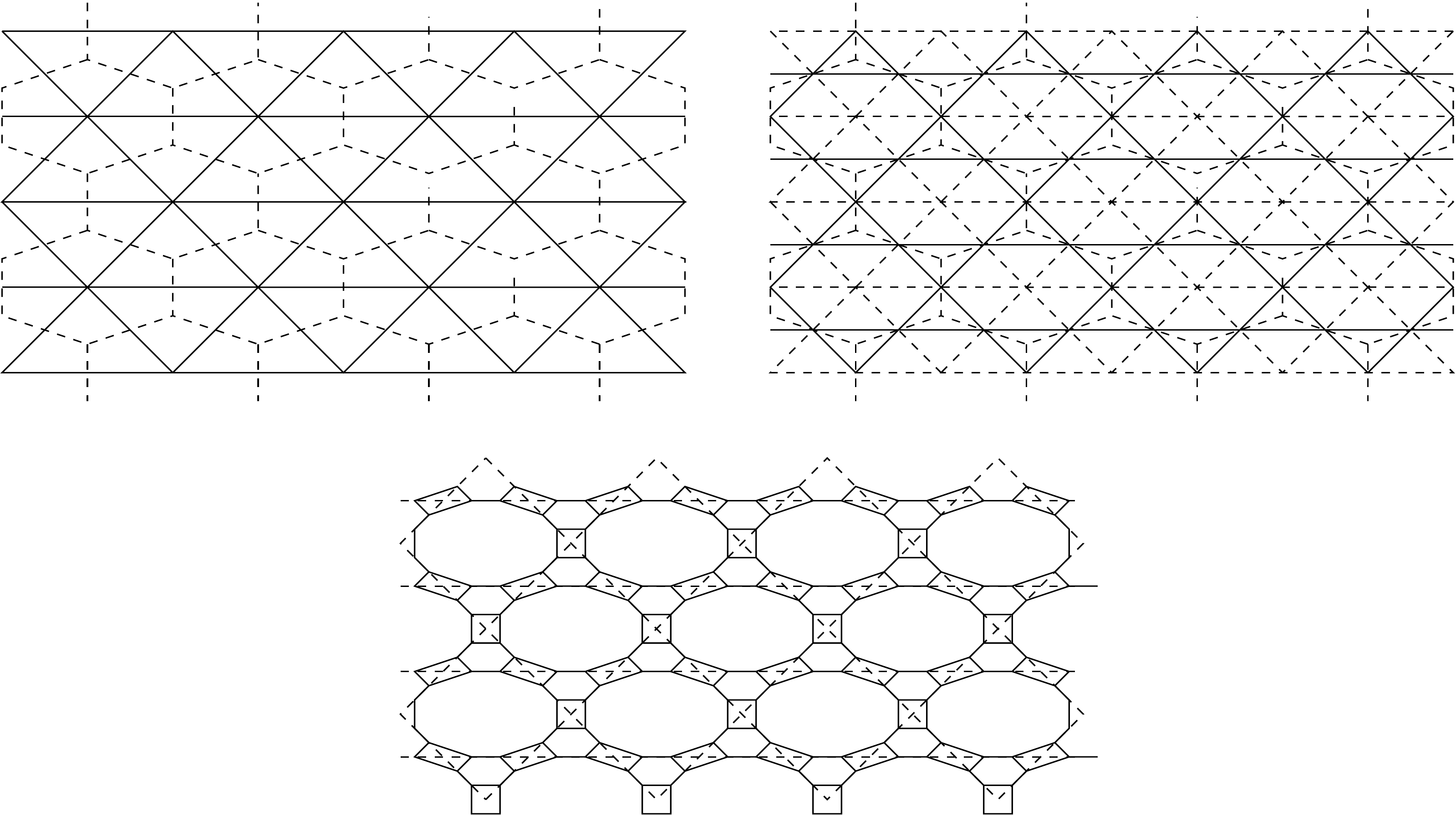}}
\end{center}
\caption{Triangular lattice $\Gamma$ (solid) and its dual hexagonal lattice (dashed); associated 4-regular graph $M$ (Kagom\'e lattice); associated bipartite graph $C$ (4-6-12 archimedean lattice)
}
%
%
\label{fig:archi}
\end{figure}

It is obvious that the correspondence is measure-preserving if the 6V weights are given by
$$(\omega_{12},\omega_{34},\omega_{56})=(\sin\theta,\cos\theta,1)$$
(where $\theta$ may depend on the city; the weight of the corresponding edge is $w=\tan(\theta/2)$) and the dimer weights are 1 for each road, $\sin\theta$ for horizontal city streets and $\cos\theta$ for vertical city streets (in the coordinates of Figure \ref{fig:6Vcity}). The various graphs derived from $\Gamma$ are illustrated in Figure \ref{fig:medial} in the generic case and in Figure \ref{fig:archi} when $\Gamma$ is the triangular lattice.

Associated to the 6V model or its bipartite dimer representation is a {\em height function} $\phi$ which may be described as follows. To each face $f$ of $M$ one associates a height $\phi(f)$ such that $\phi(f')-\phi(f)=\pi $ (resp. $-\pi$) is crossed from left to right by an oriented edge in the 6V configuration. The 6V rule ensures that $\sum_{i=0}^3(f_{i+1}-f_i)=0$, where $f_0,\dots,f_3$ are the cyclically ordered faces around a vertex $v\in M$. This gives a consistent definition of $\phi$ (up to an additive constant) on a fundamental domain. On the torus, $\phi$ is additively multivalued.

In terms of the height function, there are two natural (at least in the present context) types of insertions: electric and magnetic charges. If $f$ is a face of $M$, the r.v. $e^{i\alpha\phi(f)}$ is referred to as an electric charge $\alpha$ inserted at $f$. A pair of opposite charges $e^{i\alpha\phi(f')}e^{-i\alpha\phi(f)}$ does not depend on the choice of a base value for $\phi$. We may choose the additive constant of $\phi$ so that $\phi(f)=0\mod 2\pi$ if $f$ is a black face of $M$ (corresponding to a vertex of $\Gamma$); then $\phi(f)=\pi\mod 2\pi$ if $f$ is a white face (corresponding to a vertex of $\Gamma^\dg$). If $f\in\Gamma$, 
$$e^{\frac i2\phi(f)}=e^{-\frac i2\phi(f)}=\cos(\phi(f)/2)$$
and if $f\in\Gamma^\dg$,
$$-ie^{\frac i2\phi(f)}=ie^{-\frac i2\phi(f)}=\sin(\phi(f)/2)$$
An order line is a product $\prod_{e\in\gamma}\nu(e^\dg)$, where $\gamma$ is a path from $f$ to $f'$ on $M^\dg$. We have:
$$\sigma(f)\sigma(f')=\prod_{e\in\gamma}\nu(e^\dg)=e^{\frac i2(\phi(f')-\phi(f))}$$
provided $f,f'$ are of the same color (in the bipartite coloring of $M^\dg$). Finally we get the identification: $\sigma(f)=\cos(\phi(f)/2)$ if $f\in\Gamma$ and $\sigma(f')=\sin(\phi(f)/2)$ if $f\in\Gamma^\dg$.

Inserting a disorder variable at $e\in E_M$ amounts to introducing a {\em magnetic charge} $\pm 1$: the height function becomes locally multivalued, picking an additive constant $\pm \pi$ (depending on the type of defect: source or sink) when cycling counterclockwise around the defect. Remark that this has always to be compensated by another  (opposite) magnetic charge elsewhere.

In terms of the dimer mapping, disorder variables may be interpreted as {\em monomers} (\cite{FisSte2}). In the local correspondence, a 6V configuration with a sink defect (magnetic charge -1) at $e\in E_M\simeq (bw)\in E_C$ corresponds to a modified dimer configuration in which $b$ is matched to $e$ (the midpoint of $(bw)$) while $w$ is matched in its city. Conversely, a source defect corresponds to $w$ being matched to $e$ and $b$ matched in its city. Removing the half-dimer, one gets the more classical formulation of a monomer defect at $b$ or $w$. 

Specifically, if $\langle.\rangle_{\rm dimer}$ denotes the unnormalized expectation (ie $\langle 1\rangle_{\rm dimer}$ is the partition function) under the dimer measure on $C$, we denote
$$\left\langle.\prod_{i=1}^n{\mc O}_1(b_i){\mc O}_{-1}(w_i)\right\rangle_{\rm dimer}
$$
the unnormalized dimer measure on $C\setminus\{b_1,w_1,\dots,b_n,w_n\}$, where $b_1,\dots,b_n$ (resp. $w_1,\dots,w_n$) are $n$ black (resp. white) vertices on $C$ (using a magnetic operator notation, where as usual disorder variables represent a modification of the state space). If the number of black and white monomers does not match, we set $\prod_{i=1}^n{\mc O}_1(b_i)\prod_{j=1}^m{\mc O}_{-1}(w_j)=0$ (as the state space is then empty).

With these conventions, we have $\xi(e)={\mc O}_1(b)+{\mc O}_{-1}(w)$ and $\nu(e^-)\xi(e)=-\nu(e^+)\xi(e)={\mc O}_1(b)-{\mc O}_1(w)$.
where $e^-$ is the half-edge of $e$ containing $b$.

{\bf Correlators.}

We are now concerned with tracking correlators of order and disorder variables through the mappings (Ising, dual Ising)$\rightarrow$ 8V
$\rightarrow$ 8V$\rightarrow$ 6V$\rightarrow$ dimers.

In the Ising model on $\Gamma=(V,E)$, the basic order variable is a spin variable $\sigma_v$. By Kramers-Wannier duality, such order variables correspond to disorder ``variables" in the dual Ising model. Following Kadanoff and Ceva, a pair of disorder variables $\mu_f\mu_{f'}$ is represented by a simple path $\gamma$ from $f$ to $f'$ on $\Gamma^\dg$ (a disorder line); the weight of a configuration is then modified to
$$w'((\sigma_v)_{v\in V})=\prod_{e=(vv')\in E, e^\dg\notin\gamma,\sigma_v\sigma_{v'}=-1}w(e)\prod_{e=(vv')\in E, e^\dg\in\gamma,\sigma_v\sigma_{v'}=1}w(e)
$$
A general order-disorder correlator for the Ising model on $\Gamma$ is written as:
$$\left\langle \prod_{i=1}^{2n}\sigma(v_i)\prod_{j=1}^{2m}\mu(f_j)\right\rangle_{\Gamma}^{hv}$$
where $hv\in\{p,a\}^2$ designates the periodicity type of the spin configuration, and the $f_j$'s are connected pairwise by disjoint disorder lines drawn in a fixed fundamental domain (associated to the choice $A,B$ of homology basis). We also assume that all insertions $v_i$, $f_j$ are disjoint (at no cost since $\sigma_v^2=1$ and $\mu_f^2=1$).
 
Let us consider simultaneously a dual Ising model on $\Gamma^\dg$, with a general correlator of the form:
$$\left\langle \prod_{i=1}^{2n'}\sigma^\dg(f'_i)\prod_{j=1}^{2m'}\mu^\dg(v'_j)\right\rangle_{\Gamma^\dg}^{hv}$$
For definiteness, let us assume that the spin sites $v_i$ (resp. $f'_i$) are paired by simple ``order lines" on $\Gamma$ (resp. $\Gamma^\dg$), and that for each of the two configurations, order and disorder are disjoint (and are also disjoint from $A,B$ cycles bounding a fundamental domain).

It is clear that in the mapping (Ising, dual Ising)$\rightarrow$ 8V, an order line from, say, $v_1$ to $v_2$ on $\Gamma$ can be represented by an order line, ie a simple path $\gamma'$ from $v_1$ to $v_2$ on $M^\dg$. This is simply saying that 
$$\sigma(v_1)\sigma(v_2)=\prod_{e\in\gamma'}\nu(e^\dg)$$
Let us now consider an oriented disorder line $\gamma$ from $f_1$ to $f_2$ on $\Gamma^\dg$, say. One way to represent it goes as follows: let $e$ be an oriented edge on $\gamma$, which separates the site $v$ on the left from the site $v'$ on the right. We add a new site $(ve)$ on $e^\dg$, close to $v$, and set a spin variable $\sigma(ve)=-\sigma(v)$. Then the factor of the configuration weight corresponding to $(vv')$ is $e^{\beta J_e\sigma(ve)\sigma(v')}$ (rather than $e^{\beta J_e\sigma(v)\sigma(v')}$ in the absence of disorder).

\begin{figure}[htb]
\begin{center}
\leavevmode
\includegraphics[width=0.8\textwidth]{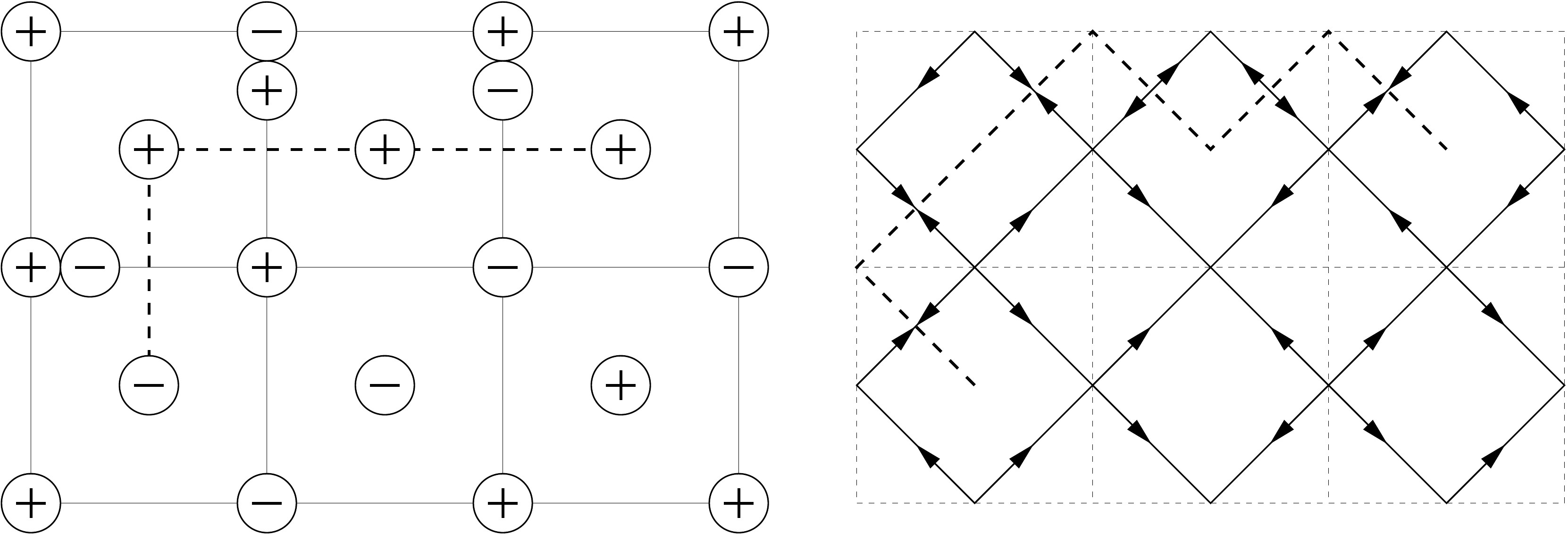}
\end{center}
\caption{Left panel: square lattice (solid), primal and dual spin configurations, disorder line (dashed, oriented from bottom left to top right). Right panel: corresponding 8V configuration with disorder line
(dashed)}
\label{fig:8Vdisorder}
\end{figure}

Correspondingly, if $e$ is an edge on $M$ which ends on $\gamma$ and is on its LHS (corresponding to $v$ in $V$ and $f$ in $F$), we define two opposite edge spin variables: $\sigma(v)\sigma^\dg(f)$ and $\sigma(ve)\sigma^\dg(f)$; in other words, we insert a disorder variable $\xi(e)$. Consequently, we can represent the Ising disorder line $\gamma$ on $\Gamma^\dg$ as an 8V disorder line $\tilde\gamma$ on $M^\dg$, which tracks $\gamma$ on its LHS. See Figure 
\ref{fig:8Vdisorder} for an illustration of the situation.

Plainly, we can recover the spin configuration on $M^\dg$ from the 8V configuration up to global spin flip. Thus we denote by $\sigma(f)$ (resp. $\mu(f)$) the endpoint of an 8V order (resp. disorder) line, where $f$ is a face of $M$, corresponding to a vertex of $\Gamma$ or $\Gamma^\dg$. For given endpoints, changing the pairings of insertions or moving the disorder lines may only change the sign of correlators.

Hence we obtain the first identity:
\begin{equation}\label{eq:boson1}
\begin{split}
\left\langle \prod_{i=1}^{2n}\sigma(v_i)\prod_{j=1}^{2m}\mu(f_j)\right\rangle_{\Gamma}^{hv}\left\langle \prod_{i=1}^{2n'}\sigma^\dg(f'_i)\prod_{j=1}^{2m'}\mu^\dg(v'_j)\right\rangle_{\Gamma^\dg}^{hv}&=2\left\langle\ind_{\eps_A=h,\eps_B=v}
\prod_{i=1}^{n+n'}\prod_{e\in\gamma_i}\nu(e^\dg)\prod_{j=1}^{m+m'}\prod_{e\in\gamma'_j}\xi(e^\dg)\right\rangle_{8V}\\
&=2\left\langle\ind_{\eps_A=h,\eps_B=v}
 \prod_{i=1}^{2n}\sigma(v_i)\prod_{i=1}^{2n'}\sigma(f'_i)
\prod_{j=1}^{2m}\mu(f_j)\prod_{j=1}^{2m'}\mu(v'_j)
\right\rangle_{8V}
\end{split}
\end{equation}
Here the 8V weights are as in \eqref{eq:8Vweights}; the $\gamma_i$'s (resp. $\gamma_j'$) are order (resp. disorder) lines on $M^\dg$ pairing the order (resp. disorder) insertions $v_1,v_{2n},f'_1,f'_{2n'}$ (resp. $f_1,\dots,f_{2m},v'_1,\dots,v'_{2m'}$). The signs $\eps_A,\eps_B$ are the products of $8V$ edge variables $\nu$ along $A$, $B$ cycles on $M^\dg$. We identify $p=1$, $a=-1$ for periodic and antiperiodic boundary conditions. The factor 2 accounts for the spin flip symmetry $(\sigma,\sigma^\dg)\leftrightarrow(-\sigma,-\sigma^\dg)$. 


The next step is to change the 8V weights from \eqref{eq:8Vweights} to \eqref{eq:8Vweights2}. This is somewhat complicated by the presence of disorder lines. One can decompose an 8V configuration with disorders (as in Figure \ref{fig:8Vdisorder}) as a collection of oriented curves: four curves are created at a source (type 8), four are destroyed at a sink (type 7); at a disorder insertion, two curves are created or two at destroyed. Denoting by $p_+$, $p_-$ the number of edge sources and sinks, $N_i$ the number of vertices of type $i$, we have:
$$4N_8+2p_+=4N_7+2p_-$$
Note that $p_+=p_-\mod 2$, since all our disorder lines on $M^\dg$ (ie those coming from Ising disorder lines on $\Gamma$ or $\Gamma^\dg$) have even length. Hence
$$\langle\cdot\rangle_{8V'}=\langle\cdot (-1)^{\frac{p_+-p_-}2}\rangle_{8V}$$
where 8V weights on the LHS (resp. RHS) are given by \eqref{eq:8Vweights2} (resp. \eqref{eq:8Vweights}). It is easily checked that if $\gamma$ is a disorder line, $p_+(\gamma),p_-(\gamma)$ the number of edge sources/sinks along this line, and $\tilde\gamma$ is an (edge disjoint) order line with the same endpoints, then:
$$(-1)^{\frac{p_+(\gamma)-p_-(\gamma)}2}=\prod_{e\in\tilde\gamma}\nu(e^\dg)$$
Consequently:
\begin{equation}\label{eq:boson2}
\begin{split}
\left\langle \prod_{i=1}^{2n}\sigma(v_i)\prod_{j=1}^{2m}\mu(f_j)\right\rangle_{\Gamma}^{hv}\left\langle \prod_{i=1}^{2n'}\sigma^\dg(f'_i)\prod_{j=1}^{2m'}\mu^\dg(v'_j)\right\rangle_{\Gamma^\dg}^{hv}
&=2\left\langle\ind_{\eps_A=h,\eps_B=v}
\prod_{i=1}^{n+n'}\prod_{e\in\gamma_i}\nu(e^\dg)\prod_{j=1}^{m+m'}\prod_{e\in\gamma'_j}\xi(e^\dg)\prod_{e\in\tilde\gamma_j}\nu(e^\dg)
\right\rangle_{8V'}\\
&\hspace{-1cm}=2\left\langle\ind_{\eps_A=h,\eps_B=v}
 \prod_{i=1}^{2n}\sigma(v_i)\prod_{i=1}^{2n'}\sigma(f'_i)
\prod_{j=1}^{2m}\mu(f_j)\sigma(f_j)\prod_{j=1}^{2m'}\mu(v'_j)\sigma(v'_j)
\right\rangle_{8V'}
\end{split}
\end{equation}
where $\tilde\gamma_j$ is an order line with the same endpoints as $\gamma'_j$, and the 8V weights are as in \eqref{eq:8Vweights2}.

Now the 8V duality simply exchanges order and disorder insertions, so that:
\begin{equation}\label{eq:boson3}
\begin{split}
\sum_{h,v=\pm1}\left\langle \prod_{i=1}^{2n}\sigma(v_i)\prod_{j=1}^{2m}\mu(f_j)\right\rangle_{\Gamma}^{hv}\left\langle \prod_{i=1}^{2n'}\sigma^\dg(f'_i)\prod_{j=1}^{2m'}\mu^\dg(v'_j)\right\rangle_{\Gamma^\dg}^{hv}&=2\left\langle
\prod_{i=1}^{n+n'}\prod_{e\in\gamma_i}\xi(e^\dg)\prod_{j=1}^{m+m'}\prod_{e\in\gamma'_j}\nu(e^\dg)\prod_{e\in\tilde\gamma_j}\xi(e^\dg)
\right\rangle_{6V}\\
&\hspace{-1cm}=2\left\langle
 \prod_{i=1}^{2n}\mu(v_i)\prod_{i=1}^{2n'}\mu(f'_i)
\prod_{j=1}^{2m}\sigma(f_j)\mu(f_j)\prod_{j=1}^{2m'}\sigma(v'_j)\mu(v'_j)
\right\rangle_{6V}
\end{split}
\end{equation}
where the 6V weights are given by \eqref{eq:6Vweights}.

\subsection{Spin structures}

We would like to obtain a bosonization identity of type \eqref{eq:boson3} which a factorized LHS. For this purpose, we need to discuss in greater details the relation between boundary conditions in the various models considered.

First we carry out the discussion in the absence of insertions. To an 8V (a fortiori 6V) configuration, we can associate a pair of signs $(\eps_A,\eps_B)$, where $\eps_{[\gamma]}=\prod_{e\in\gamma}\nu(e^{\dg})\in\{\pm 1\}$ where $\gamma$ is a closed cycle on $M^\dg$ and $[\gamma]$ is its homology class in $H_1(\Sigma,\Z/2\Z)\simeq H_1(\Sigma,\{\pm 1\})\simeq\{1,-1\}^2$. By the 8V condition, the RHS depends on $\gamma$ only through $[\gamma]$, and we have $\eps_{[\gamma_1]+[\gamma_2]}=\eps_{[\gamma_1]}\eps_{[\gamma_2]}$. 

We have 
$$\langle 1\rangle_{6V}=\langle 1\rangle_{8V}=\frac 12\sum_{h,v=\pm 1}\langle 1\rangle_{\Gamma}^{hv}\langle 1\rangle_{\Gamma^\dg}^{hv}$$
and would like to suitably twist this identity by a character of $H_1(\Sigma,\Z/2\Z)$ (a {\em spin structure}). If $\gamma,\tilde\gamma$ are two homotopic simple cycles on $M^\dg$, we have:
$$\left\langle \prod_{e\in\gamma}\nu(e^\dg)\prod_{e\in\tilde\gamma}\xi(e^\dg)\right\rangle_{8V}=\left\langle \prod_{e\in\tilde\gamma}\xi(e^\dg)\right\rangle_{8V'}=\left\langle \prod_{e\in\tilde\gamma}\nu(e^\dg)\right\rangle_{6V}=\left\langle\eps_{[\gamma]}\right\rangle_{6V}$$
For a spin configuration $\sigma$ (with periodicity conditions in $\{\pm 1\}^2\simeq\{p,a\}^2$) on $\Gamma$ or $\Gamma^\dg$, we set 
$\eps_{[\gamma]}(\sigma)=\eps_{[\gamma]}(P)$ where $P$ is the  corresponding low-temperature polygon; this depends solely on the boundary condition. This gives an interpretation of (closed) order cycles: $\eps_{[\gamma]}(\sigma)=\prod_{e=(vv')\in\gamma}\sigma(v')^{-1}\sigma(v)$ (in the presence of possibly antiperiodic boundary conditions, each factor in the RHS is unambiguously defined). If $(h,v)\in\{\pm 1\}^2$ represents the periodicity conditions, $\eps_A=h$ and $\eps_B=v$.

Dually, we can consider a disorder cycle $\gamma$ (``antiferromagnetic seams"), say drawn on $\Gamma$: the coupling for the Ising configuration on $\Gamma^\dg$ is negated on each edge crossed by $\gamma$. By drawing $\gamma$ on the boundary of a fundamental domain, introducing a disorder cycle is equivalent to shifting the periodicity condition. More precisely, let us identify $H_1(\Sigma,\{\pm 1\})$ with $\{\pm 1\}^2$ by setting $[A]=(-1,1)$, $[B]=(1,-1)$. We define an involution $[\gamma]\mapsto [\gamma]^*$ by $[A]^*=[B]$, $[B]^*=[A]$; and an ${\mb F}_2\simeq\{\pm 1\}$-bilinear pairing by $[A]\wedge [B]=[B]\wedge [A]=-1$, $[\gamma]\wedge [\gamma]=1$ for all $\gamma$. Then inserting a disorder cycle $\gamma$ shifts the periodicity type from $(h,v)$ to $(h,v)+[\gamma]^*$; and inserting an order cycle $\gamma$ introduces a sign $[\gamma]^*\wedge(h,v)$.
 
Now let us consider a spin configuration on $\Gamma$ with an order cycle $\gamma$ and a spin configuration on $\Gamma^\dg$ with a disorder cycle $\tilde\gamma$. As in \eqref{eq:boson3} (except we have closed cycles instead of open paths), we have
\begin{align*}
\langle \eps_{[\gamma]}\rangle_{6V}=
\left\langle \prod_{e\in\gamma}\nu(e^\dg)\prod_{e\in\tilde\gamma'}\xi(e^\dg)\right\rangle_{8V}=\frac 12
\sum_{h,v=\pm 1}\left\langle \eps_{[\gamma]}\right\rangle_\Gamma^{hv}
\left\langle1\right\rangle^{hv+[\gamma]^*}_{\Gamma^\dg}=\frac 12
\sum_{h,v=\pm 1}([\gamma]^*\wedge (h,v))\left\langle 1\right\rangle_\Gamma^{hv}
\left\langle1\right\rangle^{hv+[\gamma]^*}_{\Gamma^\dg}
\end{align*}
Let $q:\{\pm 1\}^2\rightarrow\{\pm 1\}$ be the quadratic form given by $q(1,1)=1$, $q(a,b)=-1$ otherwise. For $(\mu,\nu)\in\{0,1\}^2$ (a {\em sector}), we set:
$$\vphantom{\langle}^{\mu\nu}\langle\cdot\rangle_\Gamma=\sum_{h,v=\pm 1}q(h,v)h^\mu v^\nu\langle\cdot\rangle^{hv}_{\Gamma}$$
so that $\vphantom{\langle}^{\mu\nu}\langle\cdot\rangle_\Gamma$ is a signed measure with state space the disjoint union of the four spin configuration spaces corresponding to the periodicity types in $\{\pm 1\}^2$. Remark that (pairs of) order and disorder variables are defined consistently in these four spaces; and that in the low-temperature representation, there is a 2-1 correspondence between this total configuration space and the space of polygons (even degree subgraphs of $\Gamma$). Besides, we have the inversion formula
$$\langle\cdot\rangle^{hv}_{\Gamma}=\frac 14\sum_{\mu,\nu=0,1}
q(h,v)h^\mu v^\nu
(\vphantom{\langle}^{\mu\nu}\langle\cdot\rangle_\Gamma)
$$

Applying the previous identity to $[\gamma]=A,B,A+B$, we obtain:
$$\left\{
\begin{array}{ccc}
2\langle 1\rangle_{6V}  &=(00)(00)^\dg+(10)(10)^\dg+(01)(01)^\dg+(11)(11)^\dg\\
2\langle \eps_A\rangle_{6V}  &=(00)(01)^\dg-(10)(11)^\dg+(01)(00)^\dg-(11)(10)^\dg\\
2\langle \eps_B\rangle_{6V}  &=(00)(10)^\dg+(10)(00)^\dg-(01)(11)^\dg-(11)(01)^\dg\\
2\langle \eps_A\eps_B\rangle_{6V}  &=(00)(11)^\dg-(10)(01)^\dg-(01)(10)^\dg+(11)(00)^\dg
\end{array}
\right.
$$
where for legibility we write $(hv)$ (resp. $(hv)^\dg$) for $\langle 1\rangle^{hv}_{\Gamma}$ (resp. $\langle 1\rangle^{hv}_{\Gamma^\dg}$), identifying $\{\pm1\}$ with $\Z/2\Z$. Then a linear combination of these relations yields:
$$2\langle q(\eps_A,\eps_B)\rangle_{6V}=\vphantom{\langle}^{00}\langle 1\rangle_{\Gamma}\vphantom{\langle}^{00}\langle 1\rangle_{\Gamma^\dg}$$
and more generally
$$2\langle q(\eps_A,\eps_B)\eps_A^\nu\eps_B^\mu\rangle_{6V}=\vphantom{\langle}^{\mu\nu}\langle 1\rangle_{\Gamma}\vphantom{\langle}^{\mu\nu}\langle 1\rangle_{\Gamma^\dg}$$
for $\mu,\nu\in\{0,1\}$.

The argument may be carried out in the presence of insertions, as in \eqref{eq:boson3}. Indeed, as long as order and disorder lines are drawn in a fundamental domain, and the spin structures/boundary conditions are materialized by cycles on the boundary of the said domain, we have (with the $\gamma$'s as in \eqref{eq:boson3})
\begin{equation}\label{eq:boson4}
\begin{split}
{\vphantom{\left\langle \prod_{i=1}^2\right\rangle}}^{\mu\nu\!\!\!}
\left\langle \prod_{i=1}\sigma(v_i)\prod_{j=1}^{2m}\mu(f_j)\right\rangle_{\!\Gamma}
{\vphantom{\left\langle \prod_{i=1}^2\right\rangle}}^{\mu\nu\!\!\!}
\left\langle \prod_{i=1}^{2n'}\sigma^\dg(f'_i)\prod_{j=1}^{2m'}\mu^\dg(v'_j)\right\rangle_{\!\Gamma^\dg}
&=2\left\langle q(\eps_A,\eps_B)\eps_A^\nu\eps_B^\mu
\prod_{i=1}^{n+n'}\prod_{e\in\gamma_i}\xi(e^\dg)\prod_{j=1}^{m+m'}\prod_{e\in\gamma'_j}\nu(e^\dg)\prod_{e\in\tilde\gamma_j}\xi(e^\dg)
\right\rangle_{\!6V}\\
&\hspace{-1cm}=2\left\langle q(\eps_A,\eps_B)\eps_A^\nu\eps_B^\mu
\prod_{i=1}^{2n}\mu(v_i)\prod_{i=1}^{2n'}\mu(f'_i)
\prod_{j=1}^{2m}\sigma(f_j)\mu(f_j)\prod_{j=1}^{2m'}\sigma(v'_j)\mu(v'_j)
\right\rangle_{\!6V}
\end{split}
\end{equation}
for $\mu,\nu\in\{0,1\}$.

The reader will have noticed the close parallel with spin structures for dimers \cite{Kas_Ising}. In higher genus, this has been completely clarified in \cite{CimRes,CimResII}, in particular in terms of the Arf invariant, which allows to extend the above discussion to the case $g\geq 2$.

\subsection{Boundary conditions}

In simply-connected domains, we avoid difficulties related to spin structures; however a discussion of boundary conditions is needed. 

Standard boundary conditions for an Ising configuration are the wired and free conditions. A wired boundary arc is one on which spins are constant; this may be enforced by setting edge weights to zero along this arc. On a free arc, there is no interaction with the outside; for symmetry (wrt duality), this may be be represented by adding an outside vertex connected to each vertex on the free arc by an edge with weight 1. We will consider the case of a finite planar graph (thought of as the approximation of a macroscopic simply connected domain) with outer boundary partitioned into alternating wired and free boundary arcs (thought of as approximating finitely many components of the boundary).

\begin{figure}[htb]
\begin{center}
\leavevmode
\includegraphics[width=0.8\textwidth]{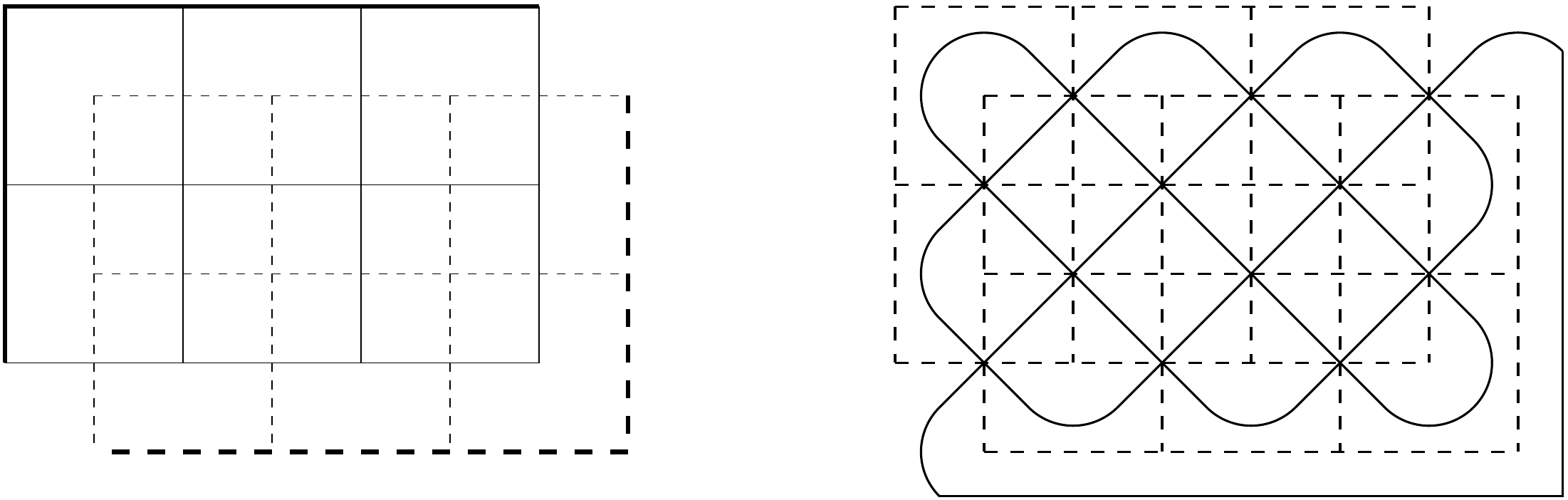}
\end{center}
\caption{Left: graph $\Gamma$ with a wired arc (bold); dual graph (dashed) with a wired arc (bold). Right: corresponding 4-regular graph $M$}
\label{fig:freewired}
\end{figure}

For simplicity, let us first consider the case where the boundary of the primal graph $\Gamma$ is partitioned in two (possibly empty) arcs which are respectively wired and free, as illustrated in Figure \ref{fig:freewired}. Correspondingly, the dual graph $\Gamma^\dg$ has a free and wired arc.

Because the spin is constant on a wired edge, we can construct the 4-regular graph $M$ (carrying corresponding the 8V model) in such a way that edges of $M$ bounce off wired edges of $\Gamma$ or $\Gamma^\dg$. Then each vertex of $M$ corresponds to a pair of dual edges with nondegenerate weights.

\begin{figure}[htb]
\begin{center}
\leavevmode
\includegraphics[width=0.8\textwidth]{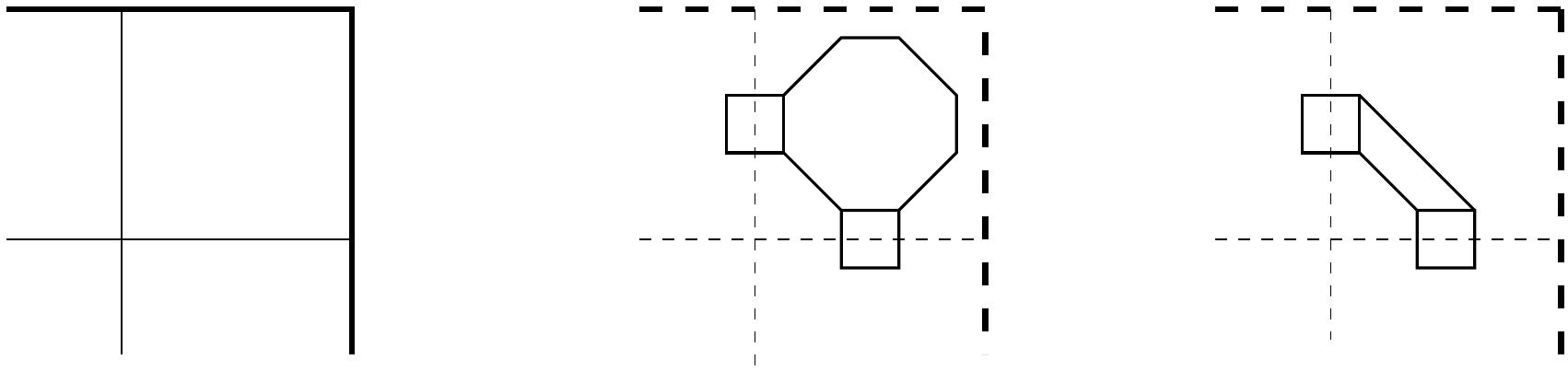}
\end{center}
\caption{Left: a face of $\Gamma$ adjacent to a wired arc (bold). Middle: corresponding graph $C$ (locally). Right: an equivalent dimer graph.}
\label{fig:wiredcorner}
\end{figure}

In the mapping to dimers, one still replaces each vertex of $V$ with a city; and these cities are connected by a road or a chain of roads (with weight 1) along each edge of $M$. Edges of $M$ now include one edge per wired edge on $\Gamma$ or $\Gamma^\dg$ (properly counted if the boundary arc is not simple), and one extended (macroscopic edge) for each pair of changes of boundary conditions. Remark that each change of boundary condition from wired to free or free to wired corresponds to a distinguished vertex of the dimer graph $C$, and that these distinguished vertices are the endpoints of the external road in $C$. This procedure produces a dimer graph $C$ which may contain chains of degree 2 vertices (in faces of $\Gamma$ or $\Gamma^\dg$ which have several consecutive edges on a wired boundary arc). The dimer model is essentially unchanged if such a chain is replaced with another chain with the same number of vertices modulo 2, see Figure \ref{fig:wiredcorner}.

The situation is somewhat more involved when the number of boundary arcs (again alternatively wired and free) is $2n>2$. Then we have the option of adding external wirings connecting some of the wired components; this external wiring materializes a given partition of the set of the $n$ wired boundary components. For instance, for two wired and two free boundary components, one can consider two cases, depending on whether or no the two wired components are connected by an external wiring. 

In order to preserve planar duality, we consider the case where the external wiring can be realized by disjoint edges, or equivalently the partition of the wired components is non-crossing. As is well known, there is a duality correspondence between the non-crossing partitions of the $n$ wired arcs and the non-crossing partitions of the dual wired arcs, and also a correspondence with the non-crossing pairings of the $2n$ endpoints of the boundary arcs. These different configurations are enumerated by Catalan's numbers. It is easy to check that a choice of non-crossing partition for the wired components on $\Gamma$ corresponds to the choice of a non-crossing pairing of the $2n$ distinguished vertices on $C$ (located at the boundary condition changes) by external extended edges. See \cite{Dub_Euler,KenWil_bound,SmiChe_Ising} for related considerations.

To summarize, given a (finite, connected) planar graph $\Gamma$ with an outer boundary consisting of $2n$ alternating wired and free arcs, and a non-crossing partition of the wired components, we can construct a 4-regular graph $M$ carrying a 6-vertex model and a bipartite graph $C$ carrying a dimer configuration in such a way that: 
\begin{equation}\label{eq:boson3plane}
\left\langle \prod_{i=1}^{2n}\sigma(v_i)\prod_{j=1}^{2m}\mu(f_j)\right\rangle_{\Gamma}\left\langle \prod_{i=1}^{2n'}\sigma^\dg(f'_i)\prod_{j=1}^{2m'}\mu^\dg(v'_j)\right\rangle_{\Gamma^\dg}=2\left\langle
\prod_{i=1}^{n+n'}\prod_{e\in\gamma_i}\xi(e^\dg)\prod_{j=1}^{m+m'}\prod_{e\in\gamma'_j}\nu(e^\dg)\prod_{e\in\tilde\gamma_j}\xi(e^\dg)
\right\rangle_{6V}
\end{equation}
where the 6V weights are given by \ref{eq:6Vweights} and order and disorder lines $\gamma,\gamma',\tilde\gamma$ are as in the toroidal case. This also applies when the whole outer boundary consists of a single wired or free arc (in which case there is no external edge in $C$).

It is also standard to consider $\pm$ boundary conditions, where the spin values on some boundary arcs are set to $+1$ or $-1$. This removes the spin flip invariance of the model, which is inconvenient for our purposes. Note however that these conditions may be realized using wired and free boundary conditions. For instance, if $\langle\cdot\rangle_{+}$ (resp. $\langle\cdot\rangle_w$) represents the unnormalized expectation for the Ising model on a planar graph $\Gamma$ with $+$ (resp. wired) boundary condition, we have
$$\langle\cdot\rangle_+=\left\langle \cdot \frac{1+\sigma(b)}2\right\rangle_w$$
where $b$ is any boundary vertex of $\Gamma$. Indeed, $\ind_{\sigma(b)=1}=\frac{1+\sigma(b)}2$. 
Note that boundary edge weights on a $\pm$ (or wired) boundary component are not counted in the spin configuration weight.

Similarly, if the boundary is split into a $+$ and a $-$ boundary arc, this may be represented in two ways: by considering a wired-free-wired-free boundary (and no external wiring between the two boundary arcs), with each free ``arc" consisting of a single edge at the endpoint of the $\pm$ arcs, and then:
$$\langle\cdot\rangle_{+,-}=\left\langle \cdot \left(\frac{1+\sigma(b_+)}2\cdot\frac {1-\sigma(b_-)}2\right)\right\rangle_w
$$
where $b_+$ (resp. $b-$) is a vertex on the $+$ (resp. $-$) boundary component. Alternatively, one can consider the boundary as a single wired component and place a pair of disorder variables at the boundary condition changes.

\subsection{Bosonization rules}

In the identities \eqref{eq:boson4},\eqref{eq:boson3plane}, we have a relation between a product of two Ising correlators involving order and disorder operators, thought of as local operators, and some 6V (or equivalently dimer) quantity, which might involve macroscopic disorder lines. This seems rather unpractical, so we will focus on the case where there is no such macroscopic line (order lines correspond to electric correlators at the endpoints). Each spin variable (on the primal graph $\Gamma$ or its dual graph $\Gamma^\dg$) corresponds the endpoint of a disorder line on the 6V model; each (primal or dual) disorder variable corresponds to the endpoint of an order line and in disorder line in the 6V model. Consequently, in these identities, if the insertions (ie order or disorder variables on the primal or dual graph) $v_i,f_j,f_i',v'_j$ are grouped in microscopic blocks of even cardinality, then the RHS is expressed in terms of local 6V (order and disorder) variables. 

In the simplest (most interesting) case, such a microscopic ``block"
 consists in $v\in\Gamma$, $f\in\Gamma^\dg$ adjacent. This corresponds to an edge $e$ in the 8V/6V graph $M$ and a road $(bw)$ in the dimer graph $C$. We denote $e^-$, $e^+$ the two half-edges of $e$, with $e^-$ corresponding to $b$. We recall that ${\mc O}_1(b)$ (resp. ${\mc O}_{-1}(w)$) represents a monomer insertion at $b$ (resp. $w$), thought of as a magnetic insertion. We also denote by $v,f$ the two corresponding faces of $M$ and by $\phi$ the 6V height field which is defined (at least as a multivalued function) on the faces of $M$. Finally, $\psi(vf)=\sigma(v)\mu(f)$ is the Kadanoff-Ceva fermion. This yields the bosonization rules recorded in Table \ref{Table:Brules}.
 
\begin{table}[htdp]
\caption{Bosonization rules}
\begin{center}
\begin{tabular}{|c|c|c|c|}
${\rm Ising}\times{\rm Ising}^\dg$,8V&8V'&6V&Dimer\\
\hline
$\sigma(v)\sigma^\dg(f)$&$\sigma(v)\sigma(f)$&
$\mu(v)\mu(f)$&
${\mc O}_1(b)+{\mc O}_{-1}(w)$\\
\hline
$\mu(f)\mu^\dg(v)$&$\mu(v)\sigma(v)\mu(f)\sigma(f)$&
$\sigma(v)\mu(v)\sigma(f)\mu(f)$&
$\pm({\mc O}_1(b)-{\mc O}_{-1}(w))$\\
\hline
$\sigma(v)\mu^\dg(v)$&$\mu(v)$&$\sigma(v)$&$\cos(\phi(v)/2)$\\
\hline
$\sigma^\dg(f)\mu(f)$&$\mu(f)$&$\sigma(f)$&$\sin(\phi(f)/2)$\\
\hline
$\psi(vf)$&$\sigma(v)\sigma(f)\mu(f)$&$\mu(v)\mu(f)\sigma(f)$&
$({\mc O}_1(b)+{\mc O}_{-1}(w))\sin(\phi(f)/2)$\\
\hline
$\psi^\dg(fv)$&$\sigma(f)\sigma(v)\mu(v)$&$\mu(f)\mu(v)\sigma(v)$&
$({\mc O}_1(b)+{\mc O}_{-1}(w))\cos(\phi(v)/2)$\\
\hline
$\psi(vf)\psi^\dg(fv)$&$\mu(f)\mu(v)$&$\sigma(f)\sigma(v)$&
$\cos(\phi(v)/2)\sin(\phi(f)/2)$\\
\hline
\end{tabular}
\end{center}
\label{Table:Brules}
\end{table}%

Here $\pm({\mc O}_1(b)-{\mc O}_{-1}(w))$ is taken to be $+({\mc O}_1(b)-{\mc O}_{-1}(w))$ at the start (say) of an 8V disorder line and $({\mc O}_1(b)-{\mc O}_{-1}(w))$ at the end (this comes from the fact that the sign of $\xi(e)\nu(e^\pm)$ depends on the half-edge $e^\pm$ on which it is evaluated).

The previous discussion may be summarized as follows. Consider dual Ising models on $\Gamma,\Gamma^\dg$, which are finite planar or toroidal graphs; to this pair is associated a weighted bipartite graph $C$. For a pair of adjacent $v\in\Gamma$, $f\in\Gamma^\dg$, we consider $X(vf)X^\dg(fv)$ a local observable of the pair of Ising configurations as in the left column of Table \ref{Table:Brules}, and $Y(bw)$ the corresponding dimer observable (in the right column). In the planar case, boundary conditions are as described earlier (which may affect the definition of $C$). By convention (for disorder variables, and spin flip symmetry for order variables), $\langle\prod_{i=1}^n X(v_if_i)\rangle_\Gamma=0$ unless the number of order and disorder variables is even. Let $c$ be the local factor:
$$c=\prod_{e\in E_\Gamma}\frac{w_e+w'_e}2$$
coming from normalizing the 6V weights \eqref{eq:6Vweights}. We have
 
\begin{Lem}[Bosonization identities]\label{Lem:Bident}
\begin{enumerate}
\item In the planar case,
$$\left\langle\prod_{i=1}^nX(v_if_i)\right\rangle_{\Gamma}\left\langle\prod_{i=1}^nX^\dg(f_iv_i)\right\rangle_{\Gamma^\dg}=2c\left\langle\prod_{i=1}
^n Y(b_iw_i)\right\rangle_{\rm dimer}$$
\item In the toroidal case, for $\mu,\nu\in\{0,1\}$,
$${\vphantom{\left\langle\prod_{i=1}^n\right\rangle}}^{\mu\nu\!\!\!\!}
\left\langle\prod_{i=1}^nX(v_if_i)\right\rangle_{\Gamma}
{\vphantom{\left\langle\prod_{i=1}^n\right\rangle}}^{\mu\nu\!\!\!\!}
\left\langle\prod_{i=1}^nX^\dg(f_iv_i)\right\rangle_{\Gamma^\dg}=2c\left\langle
q(\eps_A,\eps_B)\eps_A^\nu\eps_B^\mu
\prod_{i=1}
^n Y(b_iw_i)\right\rangle_{\rm dimer}$$
\end{enumerate}
\end{Lem}
Recall that in order to get consistent signs on both sides of these identities, we start from edge-disjoint order and disorder lines on $\Gamma,\Gamma^\dg$ and deform to obtain corresponding lines on $M^\dg$ as in Figure \ref{fig:8Vdisorder}. In some cases, in order to get local dimer observables, we replace two (or an even number of) 6V disorder lines which almost constitute a cycle by two (or the same number of) ``short" disorder lines that cross only one edge in $M$. Since the ``almost cycle" encloses an even number of Ising order and disorder variables, the rearrangement of disorder lines leaves the sign of the correlator unchanged.

This may be used in conjunction with Kramers-Wannier duality. We recall:

\begin{Lem}[Kramers-Wannier duality]
\begin{enumerate}
\item In the planar case,
$$\left\langle\prod_{i=1}^{2n}\sigma(v_i)\prod_{j=1}^{2m}\mu(f_j)\right\rangle_\Gamma=2^{|\Gamma|-1}\left(\prod_{e\in E_\Gamma}\frac{1+w_e}2\right)\left\langle\prod_{i=1}^{2n}\mu^\dg(v_i)\prod_{j=1}^{2m}\sigma^\dg(f_j)\right\rangle_{\Gamma^\dg}$$
\item In the toroidal case, if $\mu,\nu\in\{0,1\}$, 
$${\vphantom{\left\langle\prod_{i=1}^{2n}\right\rangle}}^{\mu\nu\!\!\!\!}
\left\langle\prod_{i=1}^{2n}\sigma(v_i)\prod_{j=1}^{2m}\mu(f_j)\right\rangle_\Gamma=-q(\mu,\nu)2^{|\Gamma|}\left(\prod_{e\in E_\Gamma}\frac{1+w_e}2\right)
{\vphantom{\left\langle\prod_{i=1}^{2n}\right\rangle}}^{\mu\nu\!\!\!\!}
\left\langle\prod_{i=1}^{2n}\mu^\dg(v_i)\prod_{j=1}^{2m}\sigma^\dg(f_j)\right\rangle_{\Gamma^\dg}$$
\end{enumerate}
\end{Lem}

Consequently:

\begin{Lem}[Bosonization identities: squares]\label{Lem:Bidentsq}
\begin{enumerate}
\item In the planar case,
$$\left(\left\langle\prod_{i=1}^{2n}\sigma(v_i)\prod_{j=1}^{2m}\mu(f_j)\right\rangle_{\Gamma}\right)^2=2^{|\Gamma^\dg|}\left(\prod_{e\in E_\Gamma}\frac{1+w'_e}2\right)c\left\langle\prod_{i=1}^{2n}\cos(\phi(v_i)/2)\prod_{j=1}^{2m}\sin(\phi(f_j)/2)\right\rangle_{\rm dimer}$$
\item In the toroidal case, for $\mu,\nu\in\{0,1\}$,
$$\left(
{\vphantom{\left\langle\prod_{i=1}^{2n}\right\rangle}}^{\mu\nu\!\!\!\!}
\left\langle\prod_{i=1}^{2n}\sigma(v_i)\prod_{j=1}^{2m}\mu(f_j)\right\rangle_{\Gamma}\right)^2=-q(\mu,\nu)2^{|\Gamma^\dg|}\left(\prod_{e\in E_\Gamma}\frac{1+w'_e}2\right)c\left\langle q(\eps_A,\eps_B)\eps_A^\nu\eps_B^\mu\prod_{i=1}^{2n}\cos(\phi(v_i)/2)\prod_{j=1}^{2m}\sin(\phi(f_j)/2)\right\rangle_{\rm dimer}$$
\end{enumerate}
\end{Lem}

The local prefactors are not important (as long as we ensure they do not depend on the spin structure in the toroidal case); indeed, applying the result to empty correlators (say in the planar case) yields the (more appealing):
$$\E_{\Gamma}\left(\prod_{i=1}^{2n}\sigma(v_i)\prod_{j=1}^{2m}\mu(f_j)\right)^2=\E_{\rm dimer}\left (\prod_{i=1}^{2n}\cos(\phi(v_i)/2)\prod_{j=1}^{2m}\sin(\phi(f_j)/2)\right)_{\rm dimer}$$
Other combinations may be considered, in particular the {\em energy} $\epsilon(vv')=\sigma(v)\sigma(v')$ where $v$ and $v'$ are neighboring vertices on $\Gamma$. Then in the simply connected case, if $(ff')=(vv')^\dg$, 
$$\E_\Gamma(\epsilon(vv')^2=\E_{\rm dimer}((-1)^{(\phi(v')-\phi(v))/2\pi})$$
The variable in the RHS is a sign depending on the parity of dimers on $C$ crossing $(vv')$, and may be written (somewhat suggestively) in terms of the height function as $1-(\phi(v')-\phi(v))^2/2\pi^2$. It is a local variable in the sense of dimer statistics (\cite{Ken_locstat}), which relates energy correlators (\cite{BdT_iso,BdT_per,HonSmi_energ}) to flat dimer pattern correlations (\cite{Bou_pattern}).

\subsection{The case of the square lattice}

In the (most classical) case where $\Gamma$ is modeled on the square lattice (ie is a portion of the square lattice, or the square lattice with periodic boundary conditions), the graph $C$ is the square-octogon graph. Dimers on the square-octogon graph are known to map to dimers on a square lattice by the ``urban renewal" transformation (\cite{KPW}). Because of the particular importance of the square lattice and to avoid adding a mapping to the chain, we describe directly the mapping square lattice 6V$\rightarrow$ square lattice dimer (eg \cite{FerSpo_6V}).

\begin{figure}[htb]
\begin{center}
\leavevmode
{\includegraphics[width=.8\textwidth]{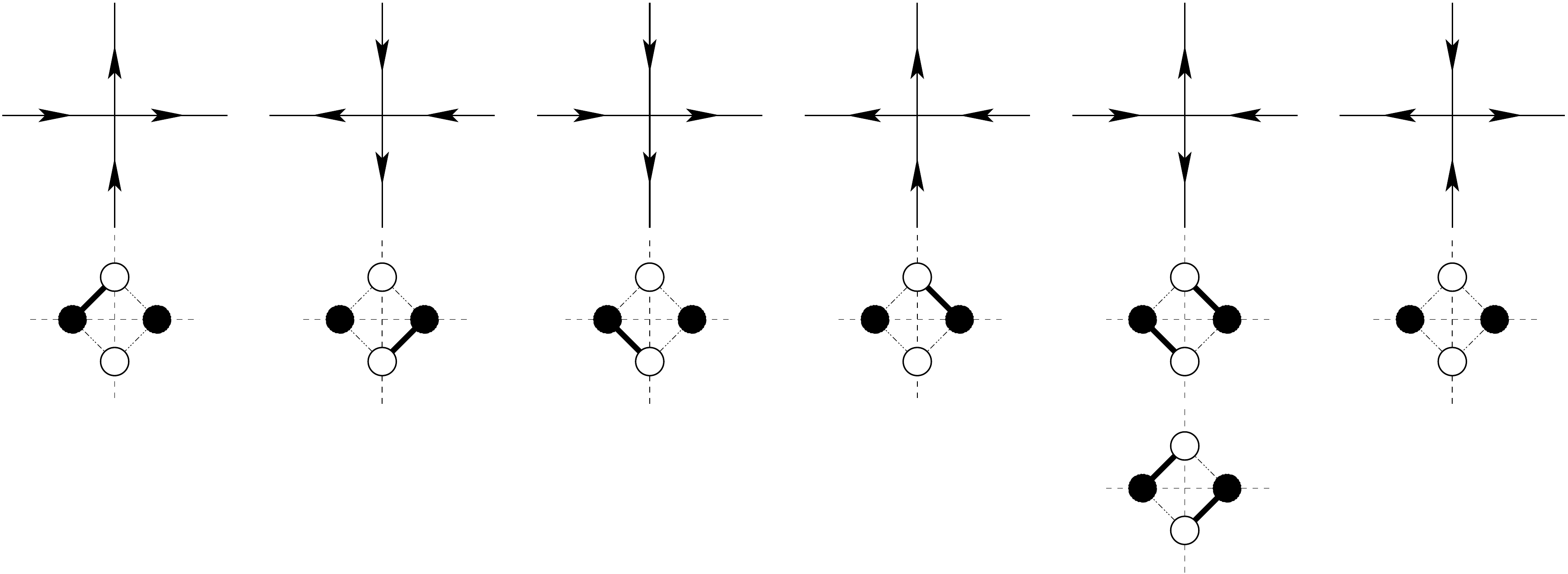}}
\end{center}
\caption{square lattice 6V configuration; associated dimer configuration
}
%
%
\label{fig:6Vdimer}
\end{figure}

Edges of $M$ are partitioned in, say, horizontal and vertical edges; their midpoints are the vertices of a square lattice $C$, with corresponding bipartite coloring. The position of the Ising square lattice $\Gamma$ relatively to the dimer square lattice is as in Figure \ref{fig:dimerorient}. The local rule is that if we orient a dimer from black to white (resp. white to black) endpoint, it makes an angle $\pm\frac\pi 4$ (resp. $\pm\frac {3\pi}4$) with the edge through $b$ (resp. $w$), with 6V orientation. See Figure \ref{fig:6Vdimer}.
Again there is a local ambiguity around type 5 vertices, and the correspondence is measure-preserving for 6V weights $(\omega_{12},\omega_{34},\omega_{56})=(\cos\theta,\sin\theta,1)$ and dimer weights $\cos(\theta)$ for NE-SW dimers and $\sin(\theta)$ for NW-SE dimers, around the generic vertex depicted in Figure \ref{fig:6Vdimer}.

In the self-dual (critical) case ($\theta=\pi/4$), we obtain the familiar uniform square lattice dimers. In the off-critical case, we have a ``flipped" weighting for dimers in the sense of \cite{Chhita_off}.

Then we have to check how electric and magnetic insertions translate in the dimer model. For electric insertions, the discussion is unchanged. For magnetic insertions, a small modification is needed. A sink defect on a horizontal 6V edge of $M$ (with midpoint $b\in C$) corresponds to a monomer defect at $b$ with charge 1 and is still denoted by ${\mc O}_1(b)$. Similarly, ${\mc O}_{-1}(w)$ denotes a source defect on a vertical edge of $M$ with midpoint $w\in C$ or equivalently a monomer at $w$. A source defect at $b$ corresponds to a {\em trimer} defect on $C$: $b$ is matched to two white vertices, one to its right and one to its left; this is denoted by ${\mc O}_{-1}(b)$, as the total charge is $-1$. The weight of the trimer is the product of the weights of the two dimers it contains. Similarly, ${\mc O}_1(w)$ denotes a trimer at $w$, corresponding to a sink defect of the 6V configuration.  

With these conventions, we record the slightly modified bosonization rules in Table \ref{Table:Brulessq}. Here $v\in\Gamma$ and $f\in\Gamma^\dg$ are adjacent; $x\in C$ is the midpoint of $[vf]$ (and may be black or white depending on the orientation of ${\vec{vf}}$). 
\begin{table}[htdp]
\caption{Bosonization rules: square lattice}
\begin{center}
\begin{tabular}{|c|c|}
${\rm Ising}\times{\rm Ising}^\dg$,8V&Dimer\\
\hline
$\sigma(v)\sigma^\dg(f)$&
${\mc O}_1(x)+{\mc O}_{-1}(x)$\\
\hline
$\mu(f)\mu^\dg(v)$&
$\pm({\mc O}_1(x)-{\mc O}_{-1}(x))$\\
\hline
$\sigma(v)\mu^\dg(v)$&$\cos(\phi(v)/2)$\\
\hline
$\sigma^\dg(f)\mu(f)$&$\sin(\phi(f)/2)$\\
\hline
$\psi(vf)$&
$({\mc O}_1(x)+{\mc O}_{-1}(x))\sin(\phi(f)/2)$\\
\hline
$\psi^\dg(fv)$&
$({\mc O}_1(x)+{\mc O}_{-1}(x))\cos(\phi(v)/2)$\\
\hline
$\psi(vf)\psi^\dg(fv)$&
$\cos(\phi(v)/2)\sin(\phi(f)/2)$\\
\hline
\end{tabular}
\end{center}
\label{Table:Brulessq}
\end{table}%
Again, the $\pm$ signs in $\pm({\mc O}_1(x)-{\mc O}_{-1}(x))$ induce a   minus per pair of such insertions (this matters only if we want to keep track of signs, which may alternatively be recovered from the asymptotic expansion).

In the plane, monomer correlators are analyzed in \cite{Dub_tors}; let us summarize these results in the case of the square lattice. If $b_1,\dots,b_n,w_1,\dots,w_n$ are $n$ black and white vertices on $\Z^2$, the monomer correlator $\Mon(b_1,\dots,b_n;w_1,\dots,w_n)$ is the limit of 
$$\frac{{\mc Z}(\Gamma_n\setminus\{b_1,\dots,w_n\})}{{\mc Z}(\Gamma_n)}=\frac{\langle\prod_i {\mc O}_1(b_i){\mc O}_{-1}(w_i)\rangle_{\rm dimer}^{\Gamma_n}}{\langle 1\rangle_{\rm dimer}^{\Gamma_n}}
$$ 
as $n$ goes to infinity, where ${\mc Z}(\Gamma)$ is the dimer partition function on the weighted graph $\Gamma$, and $\Gamma_n$ is a finite subgraph of $\Z^2$ such that $\Gamma_n\nearrow\Z^2$ and the uniform measure on dimers of $\Gamma_n$ converges as $n\rightarrow\infty$ to the maximum entropy measure on dimers on $\Z^2$. If $b_1,\dots,w_n$ have comparable pairwise distances $\gg 1$, we have
$$\Mon_{\Z^2}(b_1,\dots, w_n)\sim c_m\frac{\prod_{i\neq j}|b_i-b_j|^{1/2}|w_i-w_j|^{1/2}}{\prod_{i,j}|b_i-w_j|^{1/2}}$$
in agreement with the heuristic interpretation of monomer defects as magnetic insertions for the free field. In \cite{Dub_tors}, it is also pointed out (as a by-product of the proof) that trimer defects may also be included. If we define
$$\E_{\rm dimer}^{\Z^2}\left(\prod_{i=1}^{2n}{\mc O}_{\eps_i}(v_i)\right)=\lim_{n\rightarrow\infty}\frac{\langle\prod_{i=1}^{2n}{\mc O}_{\eps_i}(v_i)\rangle_{\rm dimer}^{\Gamma_n}}{\langle 1\rangle_{\rm dimer}^{\Gamma_n}}
$$ 
where $\eps_i=\pm 1$, $\sum_i\eps_i=0$, we have:
$$\E_{\rm dimer}^{\Z^2}\left(\prod_{i=1}^{2n}{\mc O}_{\eps_i}(v_i)\right)\sim c_m\alpha^k\prod_{i<j}|v_i-v_j|^{\eps_i\eps_j/2}$$
where $k$ is the number of trimers in the correlator, $v_i-v_j=\Theta(R)$, $R\rightarrow\infty$. In general, the constant $\alpha$ depends explicitly on the local geometry and the weighting conventions for dimers and trimers; with the the conventions in use here, $\alpha=1$. This may also be seen directly by using the 6V formulation, which is invariant under arrow reversal (which exchanges dimers and trimers).

\subsection{Isotropic Ashkin-Teller model}

Let us consider the case of the isotropic Ashkin-Teller model on $\Gamma$ (see eg \cite{Baxter_exact}). To each edge $e$ of $\Gamma$ are associated two coupling constants $J_e,J'_e$. A configuration consists in a pair of spins $\sigma(v),s(v)\in\{\pm 1\}$ for each $v\in V_\Gamma$; the weight of a configuration is
$$w(\sigma,s)=\exp\left(\beta\sum_{e=(vv')\in E_{\Gamma}} \left(J_e(\sigma_v\sigma_{v'}+s_vs_{v'})+J''_e\sigma_v\sigma_{v'}s_vs_{v'}\right)\right)$$
Notice the invariance under $(s,\sigma)\leftrightarrow (\pm s,\pm\sigma)$ and $(s,\sigma)\leftrightarrow(\sigma,s)$. Let us write
$$w_{e}(\sigma,s)=e^{\beta J_es_vs_{v'}}e^{\beta K_e\sigma_v\sigma_{v'}}$$
where $e=(vv')$, $K_e=K_e(s)=J_e+J'_es_vs_{v'}$. Following Wegner, one may apply Kramers-Wannier duality to the $\sigma$-field, regarding for now $s$ as fixed. Then
$$\sum_{\sigma:V_\Gamma\rightarrow\{\pm 1\}}\prod_ew_e(\sigma,s)=2^{|\Gamma|-1}\prod_{e=(vv')}\left(e^{\beta J_es_vs_v'}\cosh(\beta K_e)e^{-\beta K_e^\dg}\right)\sum_{t:V_{\Gamma^\dg}\rightarrow\{\pm 1\}}\prod_{e^\dg=(ff')\in E_{\Gamma^\dg}}e^{\beta K_e^\dg t_ft_{f'}}$$
where
$$\tanh(\beta K_e)=e^{-2\beta K_e^\dg}$$
If we associate an 8V configuration to the spin fields $(s,t)$ (as in Figure \ref{fig:8V}), the 8V weights are:
\begin{align*}
\omega_{12}=e^{\beta J_e}\sinh(\beta(J_e+J'_e))&&\omega_{34}=e^{-\beta J_e}\cosh(\beta(J_e-J'_e))\\
\omega_{56}=e^{\beta J_e}\cosh(\beta(J_e+J'_e))&&\omega_{78}=e^{-\beta J_e}\sinh(\beta(J_e-J'_e))
\end{align*}
Observe that:
$$\omega_{56}=\omega_{12}+\omega_{34}+\omega_{78}$$
Then replacing $\omega_{78}$ with $-\omega_{78}$ and applying $8V$ duality, we obtain the $6V$ model with weights:
\begin{align*}
\hat\omega_{12}&=\frac{\omega_{12}-\omega_{34}+\omega_{56}+\omega_{78}}4=\sinh(2\beta J_e)e^{\beta J_e'}
\\
\hat\omega_{34}&=\frac{-\omega_{12}+\omega_{34}+\omega_{56}+\omega_{78}}4=e^{-\beta J_e'}\\
\hat\omega_{56}&=\frac{\omega_{12}+\omega_{34}+\omega_{56}-\omega_{78}}4=\cosh(2\beta J_e)e^{\beta J_e'}
\end{align*}
The preceding discussion may be found in \cite{Baxter_exact}. Now this argument for partition functions may be repeated in the presence of insertions, yielding in particular the following bosonization identity:
$$\E_{AT}\left(\prod_{i=1}^{2n}\sigma(v_i)s(v_i)\right)=\E_{6V}\left(\prod_{i=1}^{2n}\cos(\phi(v_i)/2)\right)$$
for the associated 6V model with weights $(\hat\omega_{12},\hat\omega_{34},\hat\omega_{56})$. Remark that this 6V model maps to a dimer model exactly when $J'_e=0$, ie when the Askhin-Teller model decouples into two independent Ising models.

\section{Asymptotic analysis}

\subsection{Spin boundary conditions for dimers}

\begin{figure}[htb]
\begin{center}
\leavevmode
\includegraphics[width=0.8\textwidth]{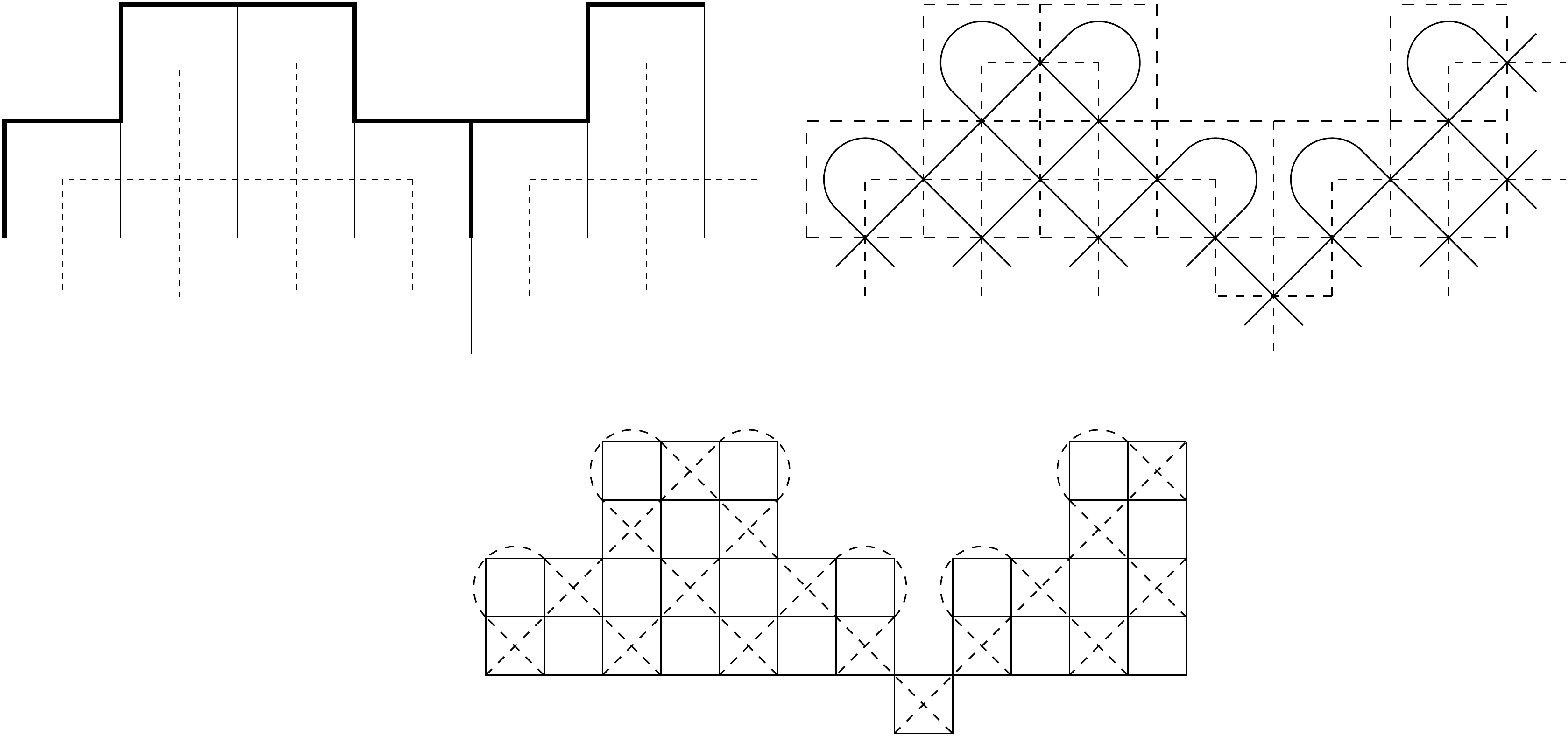}
\end{center}
\caption{Left: graph $\Gamma$ with a wired boundary (bold), dual graph (dashed). Right: corresponding 4-regular graph $M$. Bottom: dimer graph $C$.}
\label{fig:6Vdimerbound}
\end{figure}

In order to simplify the treatment of boundary conditions, we consider the square lattice case. Let us start with $\Gamma$ a simply connected portion of the square lattice, such that each face has either zero or two consecutive edges on the outer boundary, which is fully wired. The corresponding graph $M$ has two types of edges: the regular bulk edges, and edges along the boundary. In the mapping to $C$, it is rather convenient to draw three vertices on such edges, one per quarter turn. See Figure \ref{fig:6Vdimerbound}.
If we give weight 1 to the new edges (a pair of edges for each salient corner), the 6V-dimer correspondence is weight-preserving. In the isotropic case, all other edges have the same weight $\sqrt 2/2$. It is then convenient to give the same weight to these corner edges; this has the effect of multiplying the weights of all dimer configurations by a constant (an inessential {\em gauge change}).

\begin{figure}[htb]
\begin{center}
\leavevmode
\includegraphics[width=0.3\textwidth]{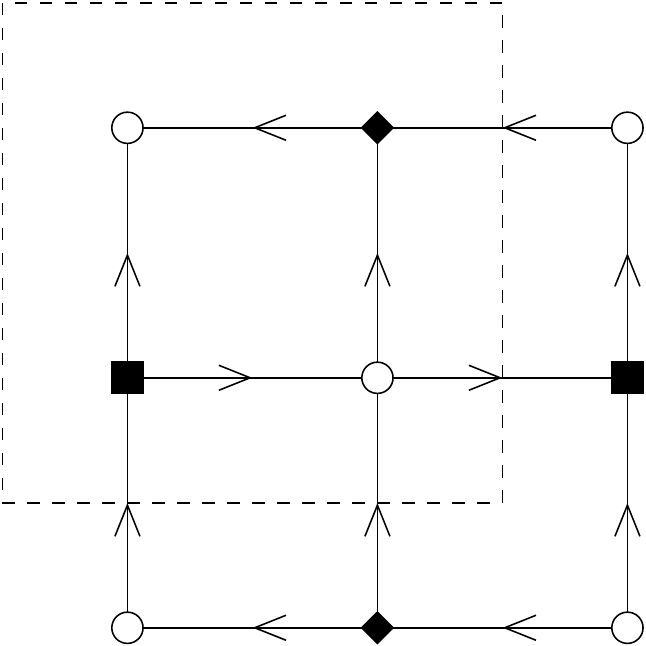}
\end{center}
\caption{A Kasteleyn orientation of the square lattice (dashed: a face of $\Gamma$)}
\label{fig:dimerorient}
\end{figure}

The graph $C$ is bipartite; we denote here $C_B,C_W$ its black and white vertices. We are concerned with the dimer model on $C$ with the special type of boundary 
In the analysis of the dimer model (\cite{Ken_IAS}), a crucial tool is the Kasteleyn operator and its inverting kernel. The Kasteleyn operator $\rK:\R^{C_B}\rightarrow\R^{C_W}$ is such that $\rK(w,b)=\pm \omega(bw)$ if $b\sim w$ ($\omega(bw)$ is the edge weight, $\sqrt 2/2$ in our normalization) and $\rK(w,b)=0$ otherwise. The sign of $\rK(w,b)$ is $1$ if $(bw)$ is oriented from $w$ to $b$ in a fixed {\em Kasteleyn orientation} of $C$ and $-1$ otherwise. A Kasteleyn orientation is such that around each face, there is an odd number of clockwise oriented boundary edges. For definiteness, we fix the orientation of Figure \ref{fig:dimerorient}, repeated periodically.
At this point let us remark that in the type of boundary considered here, the (half-decimated) height is constant on the boundary (this is easy to see from the 6V representation); by constrast, in the somewhat easier to handle Temperleyan boundary conditions, the (half-decimated) height on the boundary is proportional to the winding of the said boundary. However we retain a convenient feature of the Temperleyan case: a random walk representation of $\rK^{-1}$.

\begin{figure}[htb]
\begin{center}
\leavevmode
\includegraphics[width=0.8\textwidth]{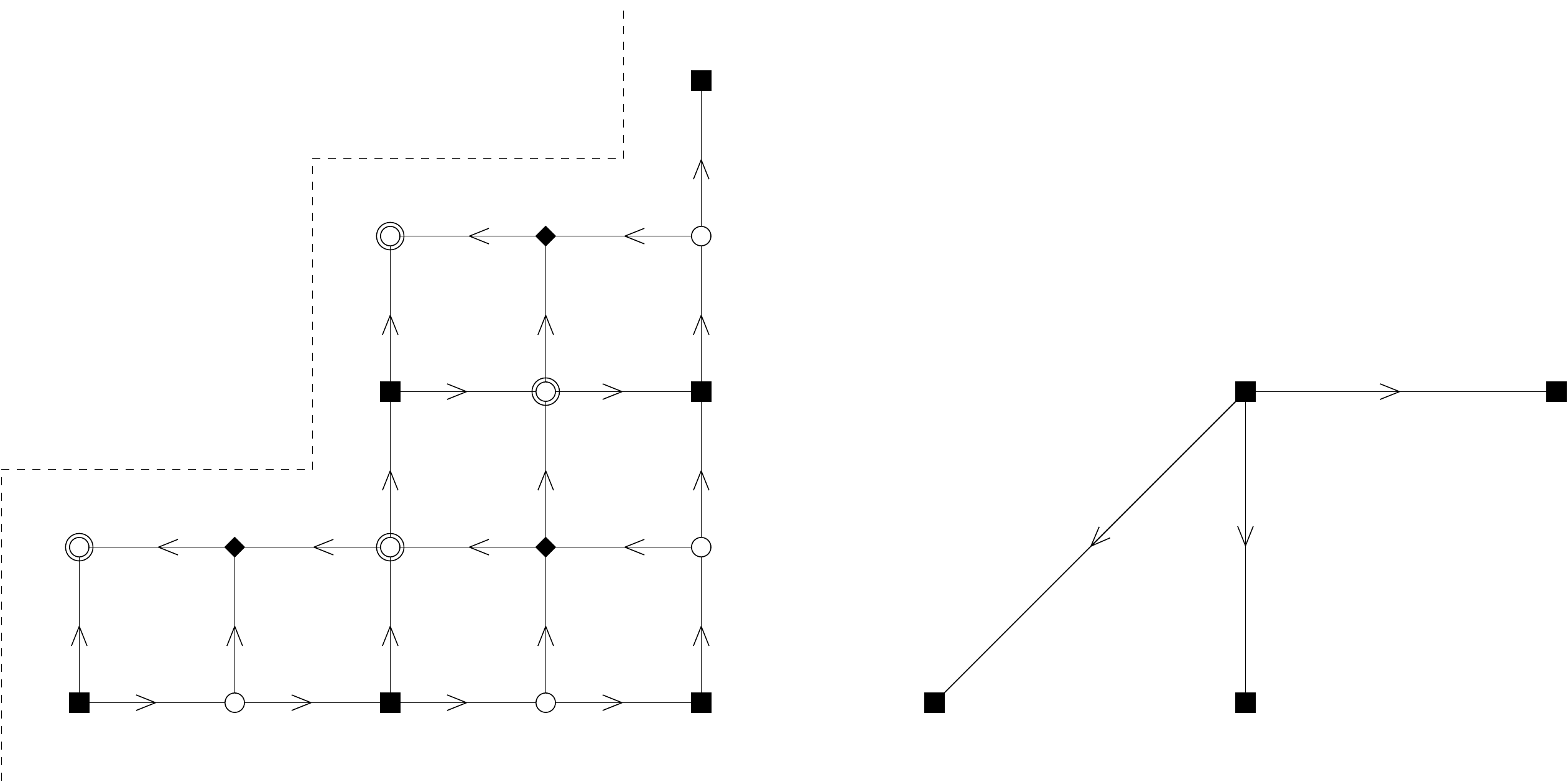}
\end{center}
\caption{Left: Oblique reflection near a SW boundary component (dashed: boundary of $\Gamma$). Right: jumps of the associated RW on $B_0$.}
\label{fig:oblique}
\end{figure}

Indeed, the black vertices may be partitioned into two subtypes, which are themselves the vertices of two dual square lattice with twice the mesh; these two types $B_0,B_1$ are represented by straight and oblique black squares in Figures \ref{fig:dimerorient}, \ref{fig:oblique}.

Fixing $w_0\in C_W$, we consider $h=\rK^{-1}(.,w_0)$ as a function on $C_{B_0}$. (By the mapping from dual Ising configurations, we know that the partition function of the dimer model is positive, and thus $\rK$ is invertible; alternatively, since the boundary height is essentially constant, Thurston's tilability criterion trivially applies). We will express it in terms of the Green kernel for a random walk on $B_0$ with suitable boundary behavior.

In the bulk, $h$ is discrete harmonic except at $b_0,b_0'$, the two points of $B_0$ adjacent to $w_0$. Indeed, $\sum_{b\sim w}\rK(w,b)h(b)=\delta_{w_0}$ by definition and $\sum_{w\sim b}\rK(w,b)(\rK f)(w)=\Lap_{B_0} f(b)$ in the bulk for any $f\in\R^{C_{B}}$ and $b\in B_0$, where $\Lap_{B_0}$ is the positive Laplacian with edge conductances $\frac 12$ on the sublattice $B_0$.

This may be adapted near the four types (pointing NE,NW,SW,SE when oriented counterclockwise) of boundary components we are considering. In all cases, a linear combination of $(\rK f)(w)$ for $w$ near $b_0$ (but not necessarily a neighbor) eliminates the dependence on $B_1$. The NW and SE pointing boundary arcs also appear in the Temperleyan case and correspond to Dirichlet (the RW is absorbed on the boundary) and Neumann (normal reflection).

On a SW boundary (Figure \ref{fig:oblique}, panel), by considering $(\rK h)(w)$ for $w$ one the the circled white vertices, one sees that $h$ is harmonic for the RW with transitions as in Figure \ref{fig:oblique}, right panel. The average jump points in the south direction. On a NE boundary, the same argument gives that $h$ is harmonic for a RW with average horizontal jump on the boundary.

Lastly, we have to consider corners. Near a salient corner, $h_{|B_0}$ and $h_{|B_1}$ are harmonic wrt to an appropriately reflected random walk on $B_0$, $B_1$. However, near a reentrant corner, $h_{|B_0}$ and $h_{|B_1}$ each fail to be harmonic at one point. Consequently, $h_{|B_0}$ (resp. $h_{|B_1}$) can be represented in terms of the Green kernel for a RW with the  appropriate reflection/absorption properties on the various boundary components, with a pair of singularities around $w_0$ and one additional singularity at the apex of each reentrant corner.

In the scaling limit, one considers a fixed polygon $P$ with boundary directions in NE,NW,SW,SE, and consider approximations of $P$ by a graph $\Gamma$ with small mesh $\delta$ and the associated graph $C$. If we consider the RW on $B_0$ with is simple in the bulk and has the transition probabilities we just described on the boundary, we get a sequence of processes which converges by standard arguments to a Brownian motion with Dirichlet/Neumann/oblique boundary condition depending on the boundary segment. Indeed, tightness is easily established and then one invokes well-posedness of the (sub)martingale problem for BM with oblique reflection in wedges (\cite{VarWil_wedge}). In the cases we consider, corners are polar (as the wedge angle is at least $\pi/2$ and the inward oblique reflection angle is at most $\pi/4$, see \cite{VarWil_wedge}). See eg \cite{Dub_RBM} for a similar reflected RW convergence argument.

It is useful to estimate the probability to exit at the apex of a reentrant corner. In the continuous limit, for an infinite wedge, the only (up to constant) nonnegative harmonic function satisfying these oblique/Neumann or oblique/Dirichlet boundary conditions can be written as $\Re(z^{-1/6})$, up to rotation and centering of the wedge. This gives estimates on exit probabilities. For instance, in the oblique/Neumann case, this extreme harmonic function is $\theta(|r|^{-1/6})$ on $C(0,r)$; it follows easily that the probability to exit the annulus $A(\lambda^{-1}r,\lambda r)$ on $C(0,\lambda^{-1}r)$ starting from $z$ with $|z|=r$ is comparable to $\lambda^{-1/6}$ for $\lambda\gg 1$. By convergence of the RW to the reflected BM, and reasoning on the crossings of $C(0,\lambda^nr_0)$, $n\in\N$, one obtains that the probability that the RW started from $z$ exits at the apex before reaching $C(0,r)$ is $\geq(|z|/r)^{-1/6-\eps}$ and $\leq (|z|/r)^{-1/6+\eps}$ for any fixed $\eps>0$ and $|z|/r$ large enough. In particular, if $h$ is discrete harmonic except at the apex and of order $1$ on the annulus $A(r,2r)$, it is $O((|z|/r)^{-1/6-\eps})$ in $D(0,r)$.

We want to estimate $\rK^{-1}_\delta(.,w_\delta)$, where $w_\delta\rightarrow w$ as the mesh $\delta\searrow 0$, with the direction of the edge of $B_0$ containing $w$ fixed. The standard argument is to establish existence and uniqueness of subsequential limits (in the proper scale). The restriction $h_\delta$ of $\rK^{-1}_\delta(.,w_\delta)$ to $B_0$ is discrete harmonic except at $b_0,b_0'$, the two vertices adjacent to $w_0$ and at one point at each reentrant corner. Fix a small $\eps>0$ and normalize $h_\delta$ so that $\sum_i||h_\delta||_{\infty, D(b_i,2\eps)\setminus D(b_i,\eps)}$ is of order 1, where the $b_i$'s are the possible singularities. From the previous argument, we have $h_\delta(z)=O(|z-c|^{-1/6-\eps})$ near a reentrant corner $c$, $h_\delta(z)=O(1/|z-w|)$ near $w$, and $h_\delta=O(1)$ elsewhere (including near salient corners), uniformly in $\delta$. By Harnack estimates, we get a bound on the Lipschitz norm of $h_\delta$ on compact subsets of $\mathring{P}$. Along a subsequential limit, $h_\delta$ converges to a harmonic function $h$ in $\mathring{P}\setminus\{w\}$. Moreover, $h$ extends continuously to the boundary arcs except possibly at corners (this may be seen from instance by the convergence of the RW to reflected BM).

In order to identify the limit, it is convenient (as in \cite{LSW_LERW}) to simultaneously consider $\hat h_\delta$, the restriction of $\rK^{-1}_\delta(.,w_\delta)$ to $B_1$. Up to extracting again a subsequence, we may also assume that $\hat h_\delta$ converges to a harmonic function $\hat h$ on $P\setminus\{w\}$, which is conjugate to $h$. Consequently, we may write $h=\psi+\bar\psi$, with $\psi$ holomorphic, and $\hat h=i\psi-i\bar\psi$. This implies that $\psi+\bar\psi=0$ on a NW segment; $\psi-\bar\psi=0$ on a SE segment; $\psi+\bar\psi=-i\bar\psi+i\psi$ on a NE segment; $\psi+\bar\psi=-i\psi+i\bar\psi$ on a SW segment. The first two express Dirichlet conditions for $h,\hat h$ respectively; the last two are obtained by examining discrete holomorphicity at a white corner vertex. 

Let $\phi:P\rightarrow\H=\{z:\Im(z)>0\}$ be a conformal equivalence; there is a single-valued determination of $\sqrt{\phi'}$ (as $\phi'$ is non-vanishing in $P$ simply connected) and the RH boundary condition can be summarized as: $\psi\parallel e^{-i\pi 8}(\phi')^{1/2}$ on the boundary. Moreover $\psi(z)=O(|z-w|^{-1})$ near $w$, $\psi$ is of order 1 near a salient corner and $\psi(z)=O(|z-c|^{-1/6-\eps})$ near a reentrant corner $c$. This essentially determines $\psi$.

Indeed, if $\psi=e^{-i\pi/8}\sqrt{\phi'}\tilde\psi$, then $\tilde\psi\circ\phi^{-1}$ is meromorphic in $\H$ with at most a simple pole at $\phi(w_0)$; is real on $\R$; has singularities of order $O(z^{-1+\eps'})$ at the images of corners, which are consequently removable; and vanishes at infinity (since $\phi'(\phi^{-1}(\infty))=0$, if we choose $\phi$ so that $\phi^{-1}(\infty)$ is an interior point of a boundary segment). Consequently,
$$\psi(z)=e^{-i\pi/8}\phi'(z)^{1/2}\left(\frac{\alpha\phi'(w)^{1/2}}{\phi(z)-\phi(w)}+\frac{\bar\alpha\bar\phi'(w)^{1/2}}{\phi(z)-\bar\phi(w)}\right)$$
where $\alpha e^{-i\pi/8}$ is the residue of $\psi$ at $w_0$. Remark that the RHS does not depend on the choice of $\phi$. The value of the residue may be identified by substracting the inverting kernel for the full plane. Thus there is a unique possible subsequential limit, which gives convergence as $\delta\searrow 0$ (for fixed $w$, uniformly in compact subsets of $P\setminus\{w\}$). Since we have convergence of $\rK^{-1}_\delta(.,w_\delta)$ for any sequence $w_\delta$ converging in $P$ (where $w_\delta$ is always on horizontal or always on vertical edges of $B_0$), this also gives uniform convergence in both variables $(w,b)$ when $w,b$ are away from the boundary and from each other.

Set $e^{i\nu(b)}=1$ (resp $i$) if $b\in B_0$ (resp. $B_1$) and $e^{i\nu(w)}=1$ (resp. $i$) if $w$ is on a horizontal (resp. vertical) edge of $B_0$. Then we have
$$\rK^{-1}_\delta(b,w)=\Re\left(\frac{e^{i\nu(b)}\phi'(z)^{1/2}}\pi\left(\frac{e^{i\nu(w)}\phi'(w)^{1/2}}{\phi(z)-\phi(w)}+\frac{e^{-i\pi/4}e^{-i\nu(w)}\bar\phi'(w)^{1/2}}{\phi(z)-\bar\phi(w)}\right)\right)+o(1)$$
uniformly in $(b,w)$ in a compact subset of $P\times P\setminus\Delta_P$.

The reader will have noticed close similarities with expressions found in \cite{Smi_ising,SmiChe_Ising,HonSmi_energ}; this is expanded on below.

\subsection{Monomers and Ising correlators}

We explain how to use the analysis of monomer correlators for the dimer model carried out in \cite{Dub_tors} to evaluate the asymptotics of Ising correlators in the planar case (for finite domains, see \cite{CHI,Dub_prep}).

Let $x_1,\dots,x_{2m},y_1,\dots,y_{2n}$ be fixed distinct points in $\C$ and $\Gamma=\delta\Z^2$. We want to estimate
$$\E_{\delta\Z^2}\left(\prod_{i=1}^{2m}\sigma(x_i^\delta)\prod_{j=1}^{2n}\mu(y^\delta_j)\right)$$
as $\delta\searrow 0$ where $\E$ is the expectation under the (unique) Gibbs measure for the Ising model on the square grid at the critical point, and $x_i^\delta-x_i=O(\delta)$ (resp. $y_i^\delta-y_i=O(\delta)$). 

For technical simplicity, let us fix a box of size $n$ (say with sides at $\pm\frac\pi 4$ angles), $n$ large enough. Let $\Gamma^n$ be the intersection of $\delta\Z^2$ with this box, say with wired boundary conditions. Let $\tilde x^\delta_i=x^\delta_i+\delta\frac {1+i}2$, and similarly $\tilde y^\delta_i=y^\delta_i+\delta\frac {1+i}2$; $u_i^\delta$ (resp. $v_i^\delta$) is the midpoint of $[x_i^\delta\tilde x_i^\delta]$ (resp.  $[x_i^\delta\tilde x_i^\delta]$). By bosonization, we have
$$\E_{\Gamma}\left(\prod_{i=1}^{2m}\sigma(x_i^\delta)\prod_{j=1}^{2n}\mu(y^\delta_j)\right)\E_{\Gamma^\dg}\left(\prod_{i=1}^{2m}\sigma(\tilde x_i^\delta)\prod_{j=1}^{2n}\mu(\tilde y^\delta_j)\right)
=(-1)^n\E_{{\rm dimer}}\left(
\prod_{i=1}^{2m}\left({\mc O}_{1}(u_i^\delta)+{\mc O}_{-1}(u_i^\delta)\right)
\prod_{j=1}^{2n}\left({\mc O}_{1}(v_i^\delta)-{\mc O}_{-1}(v_i^\delta)\right)
\right)
$$
As $n\rightarrow\infty$, the Ising measures on $\Gamma$ (with wired boundary) and $\Gamma^\dg$ (with free boundary) converge to the unique Gibbs measure for Ising on $\Z^2$ at critical temperature. Meanwhile, the dimer measure converges to the maximal entropy measure for dimers on $\Z^2$ (it follows easily from the previous section). By expanding the RHS before sending $n$ to infinity (and thus discarding contributions with non vanishing total charge), we obtain
$$\E_{\delta\Z^2}\left(\prod_{i=1}^{2m}\sigma(x_i^\delta)\prod_{j=1}^{2n}\mu(y^\delta_j)\right)^2=(-1)^n
\sum_{\stackrel{\mu_i,\nu_j=\pm 1}{\sum_i\mu_i+\sum_j\nu_j=0}}\E_{\rm dimer}\left(
\prod_{i=1}^{2m}{\mc O}_{\mu_i}(u_i^\delta)\prod_{j=1}^{2n}\nu_j{\mc O}_{\nu _j}(v_i^\delta)
\right)
$$
and consequently (\cite{Dub_tors}):

\begin{Thm}
For the Ising model on $\Z^2$ at criticality, there is $c=c_{m,n}\neq 0$ s.t. if the pairwise distances between the $x_i$, $y_j$'s are comparable and $\gg 1$, we have:
$$\E_{\Z^2}\left(\prod_{i=1}^{2m}\sigma(x_i)\prod_{j=1}^{2n}\mu(y_j)\right)^2\sim c\sum_{\stackrel{\mu_i,\nu_j=\pm 1}{\sum_i\mu_i+\sum_j\nu_j=0}}\prod_j\nu_j\prod_{i\neq i'}|x_i-x_{i'}|^{\mu_i\mu_{i'}/2}\prod_{j\neq j'}|y_j-y_{j'}|^{\nu_j\nu_{j'}/2}\prod_{i,j}|x_i-y_j|^{\mu_i\nu_j/2}
$$
\end{Thm}

Plainly, this also holds on the rectangular lattice with $Z$-invariant weights (see below). Let us also remark that, while the (pair) monomer correlation problem on $\Z^2$ (\cite{FisSte2,Hartwig_monomer}) has long been understood as closely analogous to the (pair) spin correlation problem, it is actually combinatorially equivalent to it.

\section{Relations with other approaches}

\subsection{Dimer representation}

It is well-known that an Ising configuration on $\Gamma$ can be mapped (via its high-temperature expansion) to a dimer configuration on a decoration of $\Gamma$. It is already apparent from the discussion of spin structures that there should be a direct connection between bosonization (where a pair of dual Ising configurations is in abelian duality with a bipartite dimer model) and this classical ```Pfaffian" approach (where a single Ising configuration is mapped to with a non-bipartite dimer model). We now explain this connection.

There is some flexibility in the choice of decoration: the terminal lattice introduced by Kasteleyn \cite{Kas_Ising} is the sparsest but is non-planar; in \cite{Fisher_Ising}, a planar decoration is introduced; a cyclically (around each vertex) symmetric version of the Fisher decoration is used in \cite{BdT_iso}. For our purposes, it is rather convenient to use (yet) another variant of the dimer mapping. 

\begin{figure}[htb]
\begin{center}
\leavevmode
\includegraphics[width=0.8\textwidth]{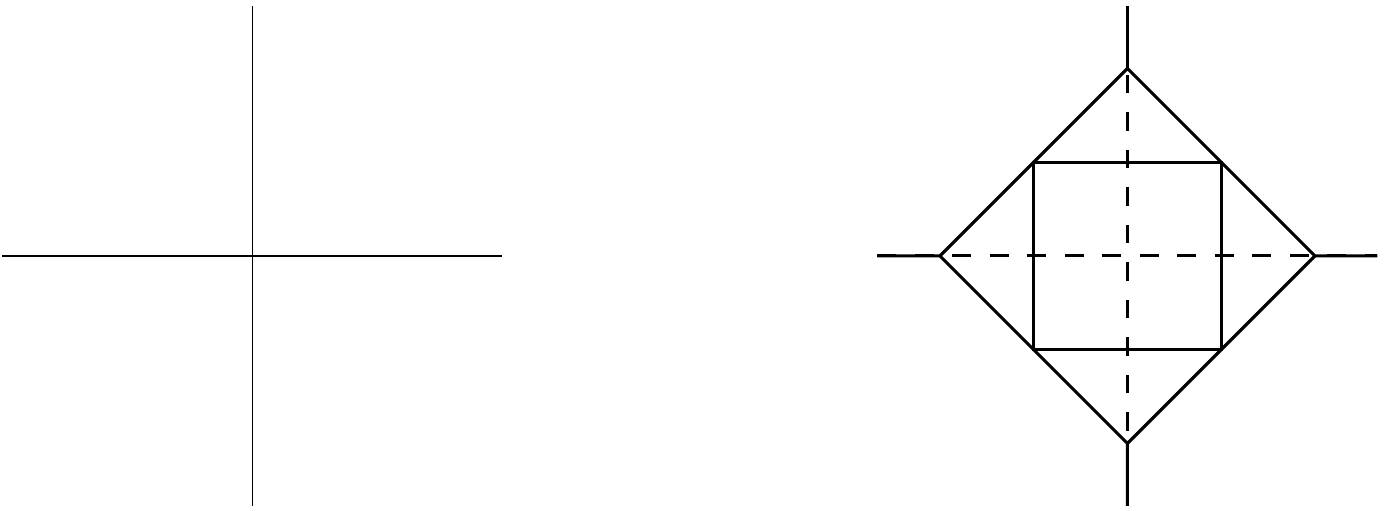}
\end{center}
\caption{A degree 4 vertex in $\Gamma$; corresponding decoration in $\Gamma_F$}
\label{fig:decor}
\end{figure}

Starting from a (planar) graph $\Gamma$, we consider a decorated graph $\Gamma_F$ obtained as follows. Each vertex $v\in\Gamma$ of degree $d$ is replaced with $2d$ vertices; one ``terminal" vertex for each outgoing edge, and one vertex for each corner of adjacent faces. Two terminal vertices are connected (by a ``road") in $\Gamma_F$ if they correspond to the same edge in $\Gamma$; two corner vertices are connected in $\Gamma_F$ if they correspond to the same vertex $v$ in $\Gamma$ and to a pair of consecutive faces around $v$; and a terminal vertex and a corner vertex are connected if the outgoing edge corresponding to the former bounds the face corner corresponding to the latter. All in all, the city replacing the vertex $v\in\Gamma$ has $2d$ vertices and $3d$ edges (streets). If $\Gamma$ is weighted, road edges of $\Gamma_F$ inherit the weights and city streets have weight 1. See Figure \ref{fig:decor}.

Let us consider a polygon, viz. an even degree subgraph $P$ of $\Gamma$. The corresponding dimer configurations ${\mf m}$ on $\Gamma_F$ are those such that each road edge is present in ${\mf m}$ iff the corresponding edge $e\in\Gamma$ is present in $P$. It is easy to check that for any given polygon $P$, there are exactly $2^{|\Gamma|}$ corresponding matchings of $\Gamma_F$. More precisely, if ${\mf m}(P)$ is the partial matching consisting of the road edges covering an edge of $P$, there are two ways to complete it to a perfect matching within each city.

A {\em Kasteleyn orientation} of a planar graph is an orientation such that around each face there is an odd number of clockwise oriented edges. The Kasteleyn operator of an oriented weighted graph $\R^{V_{\Gamma}}\rightarrow\R^{V_{\Gamma}}$ is described by its matrix elements: $K(v,v')=0$ if $v'\nsim v$; if $v\sim v'$, $K(v,v')=\omega(vv')$ if $(vv')$ is oriented from $v$ to $v'$ and $K(v,v')=-\omega(vv')$ otherwise, where $\omega(vv')$ is the weight of the (unoriented) edge $(vv')$.


We want to compare the Kasteleyn operators on $C$ and $\Gamma_F$, which are two weighted graphs constructed from the weighted graph $\Gamma$. For definiteness, we have to describe ``compatible" Kasteleyn orientations for these two graphs. Let us start with a ``geometric" Kasteleyn orientation for $C$, which is defined in terms of its embedding (this will also be convenient in the isoradial setting).

For each pair $v\in\Gamma$, $f\in\Gamma^\dg$ adjacent, choose $\nu=\nu(v,f)$ s.t. $\arg(f-v)=2\nu\mod 2\pi$ (where we somewhat abusively identify the vertex $v$ and its image in the complex plane under the graph embedding; equivalently, choose a square root of $f-v$ in $\C$). Edges in $C$ are either ``roads" or ``streets"; let us orient all roads counterclockwise around vertices of $\Gamma$. A city vertex corresponds to a pair $v\in\Gamma$, $f\in\Gamma^\dg$; a city street has endpoints of type $(vf),(v'f)$ or $(vf),(vf')$. In the first case, the street $s$ is parallel to $(vv')$ and is oriented in such a way that $(e^{i(\nu(vf)+\nu(v'f))},s)$ is direct (viz. the oriented angle is in $(0,\pi)$); in the second case, we also ask that $(e^{i(\nu(vf)+\nu(vf'))},s)$ be direct. Given the choices of $\nu$'s, this defines an orientation of $C$. Replacing $\nu(vf)$ with $\nu(vf)+\pi$ has the effect of reverting the orientation of two edges in four faces of $C$ (two city faces, and the ``large" faces corresponding to $v$ and $f$). Thus we simply need to check the Kasteleyn condition for each face of $C$ for an arbitrary local choice of $\nu$'s. This is checked directly for cities. For a large face around $v$ in $\Gamma$, let us enumerate the faces $f_1,\dots,f_d$ around $v$ in counterclockwise order and choose the $\nu(vf_i)$'s s.t. $\nu(vf_{i+1})-\nu(vf_i))\in (0,\frac\pi 2)$ for $i=1,\dots,d-1$ (this forces $\nu(vf_1)-\nu(vf_d)\in (\pi,\frac {3\pi} 2)$). Then around $v$ we have $d$ counterclockwise roads, $d-1$ ccwise streets and one clockwise street. For a large face around $f\in\Gamma^\dg$, we enumerate the vertices $v_1,\dots,v_d$ around $f$ in counterclockwise order and choose the $\nu(v_if)$'s s.t. $\nu(v_{i+1}f)-\nu(v_if))\in (0,\frac\pi 2)$ for $i=1,\dots,d-1$ (this forces $\nu(vf_1)-\nu(vf_d)\in (\pi,\frac {3\pi} 2)$). Then around $f$ there are $d$ clockwise roads, $d-1$ clockwise streets and one counterclockwise street. We refer to such a choice of $\nu$'s as standard.

We want to construct a Kasteleyn orientation for $\Gamma_F$ derived from the orientation of $C$. This may be done as follows. Each city of $\Gamma_F$ is centered at a vertex $v\in \Gamma$. If $v$ has degree $d$ (in $\Gamma$), the city $\hat v$ around $v$ has $2d$ outer streets and $d$ inner streets; the face $\tilde v$ of $C$ around $v$ has degree $2d$. There is a correspondence between outer streets of $\hat v$ and edges of $C$ around $v$; for definiteness, let us fix an edge $e\in E_\Gamma$ starting at $v$; we associate the outer street of $\hat v$ starting at the terminal vertex corresponding to $E$ and going clockwise around $v$ with the edge of $C$ around $v$ that crosses $e$; then the other $2d-1$ edges in $C$ and $\Gamma_F$ are paired cyclically.

\begin{figure}[htb]
\begin{center}
\leavevmode
\includegraphics[width=0.8\textwidth]{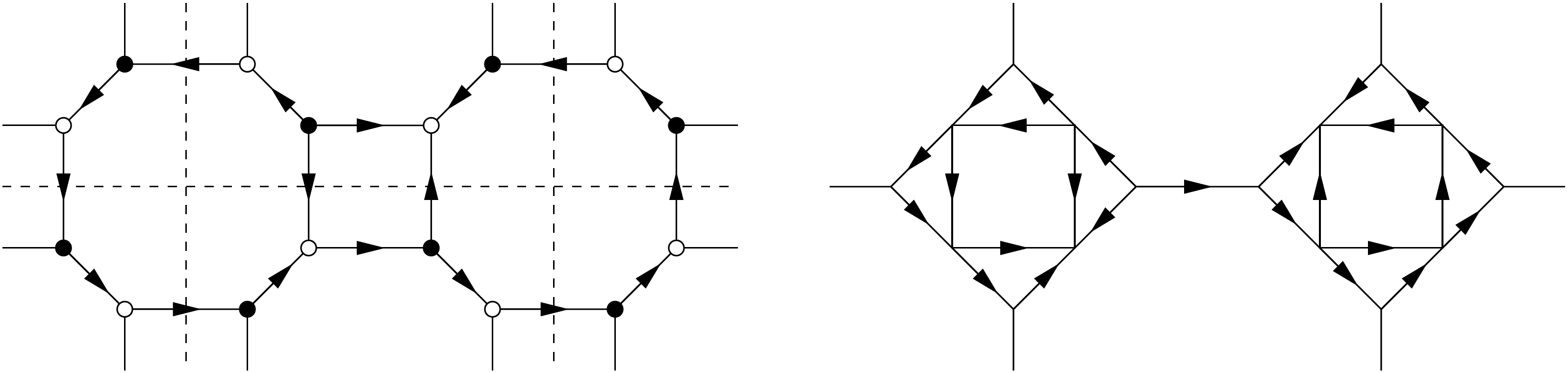}
\end{center}
\caption{Compatible standard orientations for $C$ (dashed: $\Gamma$) and $\Gamma_F$ around an edge of $\Gamma$}
\label{fig:kastorient}
\end{figure}

In this correspondence, outer streets in $\Gamma_F$ inherit an orientation from edges in $C$. Then we can orient each inner street in such a way that the Kasteleyn condition is satisfied for the outer triangle of the city containing this inner street. Since the Kasteleyn condition is satisfied for the large face of $C$ containing $v$, it is easy to see that the inner face of the city $\hat v$ in $\Gamma_F$ is clockwise odd. We still have to set the orientation of roads in $\Gamma_F$ and check the Kasteleyn condition for faces of $\Gamma_F$ which correspond to faces of $\Gamma$. Let us fix such a face $f\in\Gamma^\dg$ and $e\in E_\Gamma$ an edge on its boundary (identified with a road in $\Gamma_F$). There is a city $\tilde e$ of $C$ associated to $e$. We may assume that $e$ (oriented ccwise around $f$) goes from west to east; then the edges of $\tilde e$ are labelled W,S,E,N. If neither or both of W and S are oriented clockwise around $\tilde e$, we orient $e\in E_{\Gamma_F}$ ccwise around $f$; conversely, if exactly one of W and S is clockwise oriented, we orient $e$ clockwise around $f$. It is then immediate to check that if $f'$ is the other face of $\Gamma$ with $e$ on its boundary, we get the same orientation for $e$ reasoning around $f'$ (since the face $\tilde e$ is clockwise odd in $C$); and that the Kasteleyn condition is satisfied at the face of $\Gamma_F$ corresponding to $f$ (this uses the Kasteleyn condition for the large face of $C$ corresponding to $f$ and for the cities of $C$ corresponding to edges on the boundary of $f$). The local situation is illustrated in Figure \ref{fig:kastorient}

We want to relate the Kasteleyn operators $K_C:\R^{V_C}\rightarrow\R^{V_C}$ and $K_F:\R^{V_{\Gamma_F}}\rightarrow\R^{V_{\Gamma_F}}$. These are defined in terms of an orientation, hence in terms of a choice of $\nu(vf)$'s; changing one the $\nu$'s results in conjugating $K_C$ (resp. $K_F$) with a diagonal matrix with $\pm 1$ diagonal entries (in both cases, with exactly two $-1$ diagonal entries). Hence in local computations we may fix a ``standard" local choice of $\nu$'s as above.

Let us start from $K_F$; $\Gamma_F$ is not bipartite but its vertices are naturally partitioned evenly in two types: inner and outer vertices (seen from a vertex of $\Gamma$). An inner vertex corresponds to a pair $v\in\Gamma,f\in\Gamma^\dg$; an outer (``terminal") vertex corresponds to an endpoint of an edge $e\in\Gamma$. Let us denote by $A$ the set of inner vertices and $B$ the set of outer vertices, so that $V_{\Gamma_F}=A\sqcup B$ and $|A|=|B|$. Naturally $\R^{V_{\Gamma_F}}=\R^A\oplus\R^B$; it is also convenient to consider another decomposition of this space into, on the one hand, functions supported on $A$ and on the other hand functions $f$ which are ``holomorphic" on $B$, in the sense that $(K_F f)_{|B}=0$. Let us denote $H_B=K_F^{-1}(\R^A)\subset \R^{V_{\Gamma_F}}$ this subspace. 

Let us observe that  $\left.{K_F}\right|_{A}^B$ (restricted to $A$ and corestricted to $B$) has a diagonal block decomposition where each diagonal block (corresponding to a city) is of type
$$\left(
\begin{matrix} 
      1 & -1 & 0 & \cdots &0 \\
      0 & \ddots & \ddots &\ddots & \vdots\\
      \vdots&\ddots&\ddots&\ddots&0\\
      0 & \cdots & 0 & 1 & -1\\
      1 & 0 & \cdots & 0 & 1
   \end{matrix}
\right)$$
which has determinant 2 (irrespective of the parity of its dimension $d$, which is the degree of the corresponding vertex in $\Gamma$). (This is for a cyclic ordering of $A$ and $B$ vertices around $v$, and a standard orientation where all but one of outer streets are ccwise oriented). 

Consequently, a function $f\in H_B$ is completely determined by its restriction to $B$, and more precisely any function supported on $B$ can be extended uniquely to a function in $H_B$. In particular, if $b\in B$, we let $f_b$ be the unique element in $H_B$ s.t. $(f_b)_{|B}=\delta_b$. Plainly, the $f_b$'s constitute a basis of $H_B$, and the change of basis from the standard basis of $\R^{V_{\Gamma_F}}$ to the basis $((f_b)_{b\in B},(\delta_a)_{a\in A})$ (corresponding to the decomposition $\R^{V_{\Gamma_F}}\simeq H_B\oplus\R^A$) has determinant $\pm 1$.

Hence, if we write $K_F$ as an operator $H_B\oplus\R^A\rightarrow\R^A\oplus\R^B$ (viz. using the basis $((f_b)_{b\in B},(\delta_a)_{a\in A})$ on the LHS and the standard basis on the RHS), we get a block decomposition
$$\left(
   \begin{matrix} 
      \left.(K_F)\right|_{H_B}^A & \ast \\
      0 & \left.(K_F)\right|_{A}^B \\
   \end{matrix}
\right)$$
Since $\left.(K_F)\right|_{A}^B $ is decomposed into (explicitly) invertible local blocks, the problem of evaluating the determinant and the inverse of $K_F$ is essentially equivalent to the same problem for 
$\hat K_F\stackrel{def}{=}\left.(K_F)\right|_{H_B}^A $ . Recall that for block triangular matrices, $   \begin{pmatrix} 
      M & N \\
      0 & P \\
   \end{pmatrix}^{-1}=\begin{pmatrix} 
      M^{-1} & -M^{-1}NP^{-1} \\
      0 & P^{-1} \\
   \end{pmatrix}$.

Let us evaluate the matrix elements of $\hat K_F$. For definiteness, we work in the neighborhood of the horizontal edge $e$ of $\Gamma$ depicted in Figure \ref{fig:kastorient} (with these orientations); let $b$, $b'$ be the left and right endpoints of the corresponding road in $\Gamma_F$; let $w$ be the weight of $e$. The values of $2w^{-1}f_b$ are 
$$   \begin{matrix}
      0&\leftarrow &0&&&&1&\leftarrow &1\\
      &&&\nwarrow&&\nearrow\\
      \downarrow&&\downarrow&&2w^{-1}\longrightarrow0&&\uparrow&&\uparrow\\
      &&&\swarrow&&\searrow\\
       0&\rightarrow &0&&&&1&\rightarrow &1
   \end{matrix}$$
(where we represent $b,b'$ and the $A$ vertices of the two cities containing $b$ and $b'$; all other values are 0); similarly, the values of $-2w^{-1}f_{b'}$ are:
$$   \begin{matrix}
      1&\leftarrow &1&&&&0&\leftarrow &0\\
      &&&\nwarrow&&\nearrow\\
      \downarrow&&\downarrow&&0\longrightarrow-2w^{-1}&&\uparrow&&\uparrow\\
      &&&\swarrow&&\searrow\\
       1&\rightarrow &1&&&&0&\rightarrow &0
   \end{matrix}$$
 The extension to the general case (where the two vertices of $\Gamma$ have general degree) is straightforward. We see that $K_F f_b$ and $K_F f_{b'}$ are supported on the four A-vertices adjacent to $b$ or $b'$; let us label this vertices 
$$
 \begin{matrix}
      a_1&&&&a_4\\
      &\nwarrow&&\nearrow\\
      \downarrow&&b\longrightarrow b'&&\uparrow\\
      &\swarrow&&\searrow\\
       a_2&&&&a_3
\end{matrix}
$$
Then the restriction of $\hat K_F$ to columns corresponding to $b,b'$ and rows corresponding to $a_1,\dots,a_4$ reads:
$$   \begin{pmatrix} 
     1 & 1&-w&w \\
      w & -w & 1 & 1
   \end{pmatrix}^t
$$
Using the trigonometric parametrization $w=\tan\theta/2$, we notice the identity:
$$
 \begin{pmatrix} 
     \cos(\theta/2)&-\sin(\theta/2)\\
     \sin(\theta/2)&\cos(\theta/2)
   \end{pmatrix}
 \begin{pmatrix} 
     1 & 1&-\tan(\theta/2)&\tan(\theta/2) \\
      \tan(\theta/2) & -\tan(\theta/2) & 1 & 1
   \end{pmatrix}
=\cos(\theta/2)
\begin{pmatrix} 
     \cos(\theta)&1&-\sin(\theta)&0\\
     \sin(\theta)&0&\cos(\theta)&1
   \end{pmatrix}
$$
If we consider $\Gamma^\dg_F$ a Fisher decoration of $\Gamma^\dg$ (with dual weights \eqref{eq:KW}), there is natural 1-1 correspondence between type $A$ vertices of $\Gamma_F,\Gamma_F^\dg$ and a natural 2-2 correspondence between type $B$ vertices. At this stage it is easy to see that $\hat K_F^\dg$ may be written in terms of $\hat K_F$, yielding a version of Kramers-Wannier duality at the Fisher representation level.

Let us turn to $K_C$. As $C$ is bipartite, $\R^{V_C}\simeq\R^{V^W_C}\oplus\R^{V^B_C}$ where $V^B_C$ (resp. $V^W_C$) designates the black (resp. white) vertices of $C$. The corresponding block decomposition of $K_C$ reads   $\begin{pmatrix} 
      0 &K_C^{BW} \\
      K_C^{WB} & 0 \\
   \end{pmatrix}$
   with $K_C^{BW}:\R^{V^B_C}\rightarrow\R^{V^W_C}$ and $K_C^{WB}=-(K_C^{WB})^t$, $\Pf(K_C)=\det(K_C^{BW})$ (where $\Pf$ designates the Pfaffian of the antisymmetric matrix $K_C$ in a standard basis); this is the standard for planar bipartite graphs.
   
Here, we can establish a natural bijection between black vertices of $C$ and type $A$ vertices of  $\Gamma_F$; and between white vertices of $C$ and type $B$ vertices of $\Gamma_F$. This gives the following labels for vertices of $C$ near the edge $e$:
$$\begin{matrix}
b''&&&&&&a_4\\ 
&\nwarrow&&&&\swarrow\\
&&a_1&\rightarrow&b'\\
&&\downarrow&&\uparrow\\
&&b&\rightarrow&a_3\\
&\nearrow&&&&\searrow\\
a_2&&&&&&b'''\\ 
\end{matrix}$$
Then the submatrix of $K_C^{BW}$ with columns indexed by $a_1,\dots,a_4$ and rows indexed by $b,b'$ is
\begin{equation}\label{eq:Kholom}
\begin{pmatrix} 
     \cos(\theta)&1&-\sin(\theta)&0\\
     \sin(\theta)&0&\cos(\theta)&1
   \end{pmatrix}
\end{equation}
We conclude that there is a rotation matrix $O$ (block diagonal with $2\times 2$ blocks) and $D$ diagonal such that:
$$O\hat K_F^t=D\cdot K_C^{BW}$$
where the entries of $D$ are of type $\cos(\theta_e/2)$, with $\theta_e$ such that $w_e=\tan(\theta_e/2)$. Remark that the {\em Pfaffian} of $K_F$ enumerates (weighted) dimers on $\Gamma_F$ (corresponding to polygons on $\Gamma$), and the {\em determinant} of $K_C^{BW}$ enumerates (weighted) dimers on $C$.

Let us point out some consequences. If the graph $\Gamma$ is periodic (viz. there is a group of isomorphisms $\simeq\Z^2$ operating on $\Gamma$ with a finite quotient, say generated by a right shift $R$ and an up shift $U$), $\Gamma_F$ and $C$ are also periodic and we may choose orientations periodically. Plainly, $\hat K_F:H_B\rightarrow\R^A$ and $K^{BW}_C:\R^{V_C^B}\rightarrow\R^{V_C^W}$ and the connecting matrices $O,D$ commutes with translations. If $z,w\in\C^*$, and $V$ is a (complexified) vector space on which $R,U$ operate, let $V_{z,w}=\{v\in V: Rv=zv, Uv=wv\}$. Then $(\hat K_F)_{z,w}:(H_B)_{z,w}\rightarrow(\C^A)_{z,w}$ and $(K^{BW}_C)_{z,w}:(\C^{V_C^B})_{z,w}\rightarrow(\C^{V_C^W})_{z,w}$ are related by $O_p(\hat K_F)_{z,w}^t=D_p\cdot (K_C^{BW})_{z,w}$, where $O_p,D_p$ do not depend on $z,w$. Consequently, if we consider the characteristic polynomials $P_F(z,w)=\det((K_F)_{z,w})$ and $P_C(z,w)=\det((K_C)_{z,w})$ (this is well-defined modulo sign for a choice of standard bases for the various vector spaces), we have
$$P_F(z,w)=cP_C(z,w)$$
where $c\neq 0$ depends only on the degree distribution and the weights (the block $\left.(K_F)\right|_A^B$ operates only within cities and contributes $\pm 2^{|\Gamma_p|}$ to the constant $c$, where $\Gamma_p$ is a fundamental domain of $\Gamma$). See \cite{KOS,BdT_per,Li_Fisher} for related questions.

\begin{figure}[htb]
\begin{center}
\leavevmode
\includegraphics[width=0.6\textwidth]{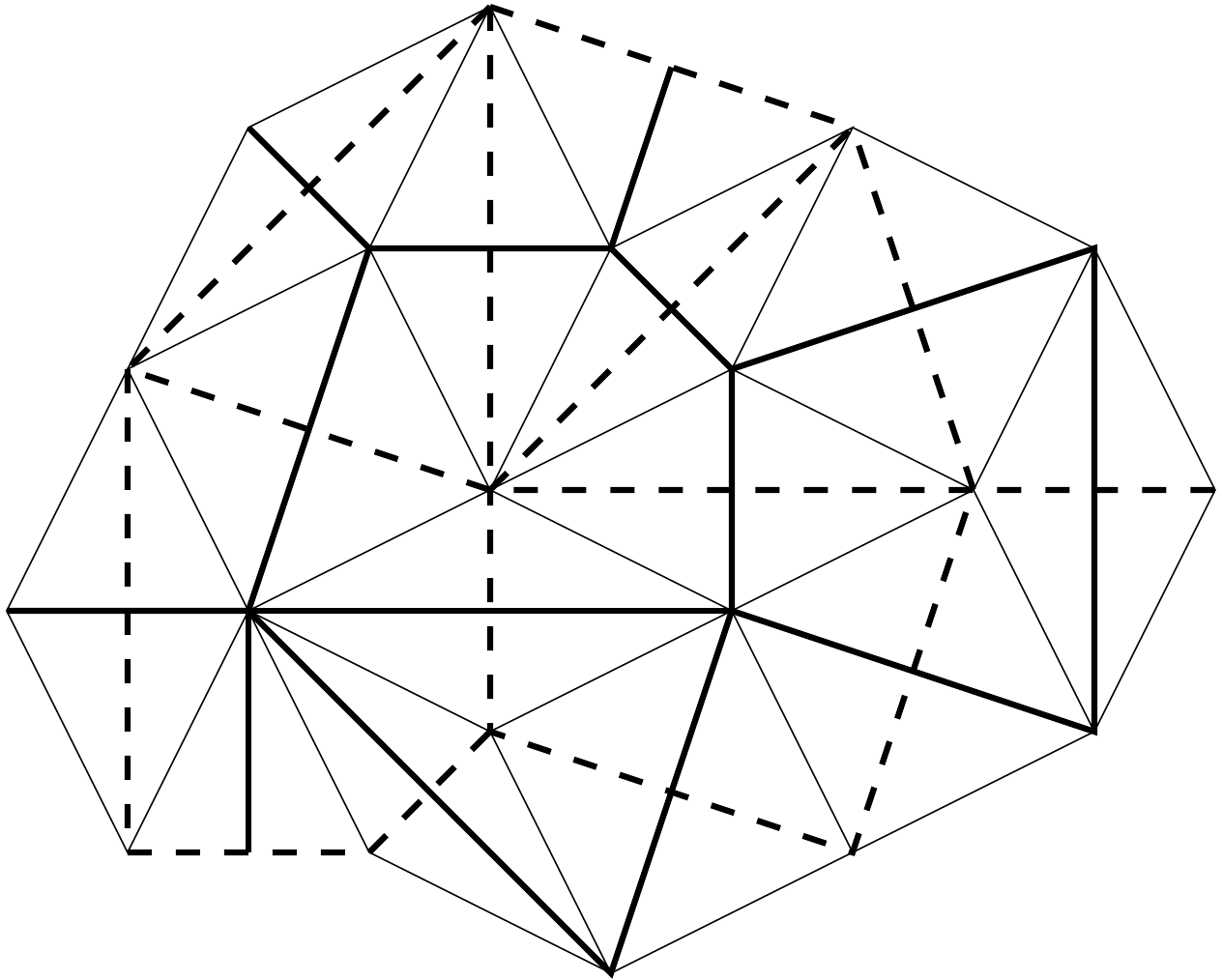}
\end{center}
\caption{a rhombus tiling $\dmd$ (thin); associated pair of dual isoradial graphs $\Gamma,\Gamma^\dg$ (bold, dashed)}
\label{fig:isoradial}
\end{figure}

A setting where $\Gamma$ is not necessarily periodic but is still tractable is the {\em isoradial} case. Start from a tiling of the plane by rhombi (the {\em diamond graph} $\dmd$), say with edge length 1 (for normalization). This is a bipartite graph. By retaining every other vertex and connecting them by a diagonal in each rhombus, one obtains an {\em isoradial} graph $\Gamma$ (its dual $\Gamma^\dg$ is also isoradial and corresponds to the other half of the rhombus tiling vertices). See Figure \ref{fig:isoradial}.
There is a special choice of Ising edge weights on $\Gamma$ associated to its isoradial embedding: the {\em critical $Z$-invariant weights} (\cite{Baxter_Zinv,BdT_iso,BdT_per,SmiChe_ising}), which are in particular significant from the Yang-Baxter point of view. If $e\in E_\Gamma$, its weight is
$$w_e=e^{-2J_e}=\tan(\theta_e/2)$$
where $\theta_e$ is half of the angle under which $e$ is seen from the center of a face of $\Gamma$ with $e$ on its boundary. Remark that on the dual graph, $w_{e^\dg}=\tan(\frac\pi 4-\frac{\theta_e}2)$, a parameterization of the Kramers-Wannier duality relation \eqref{eq:KW}.

\begin{figure}[htb]
\begin{center}
\leavevmode
\includegraphics[width=0.8\textwidth]{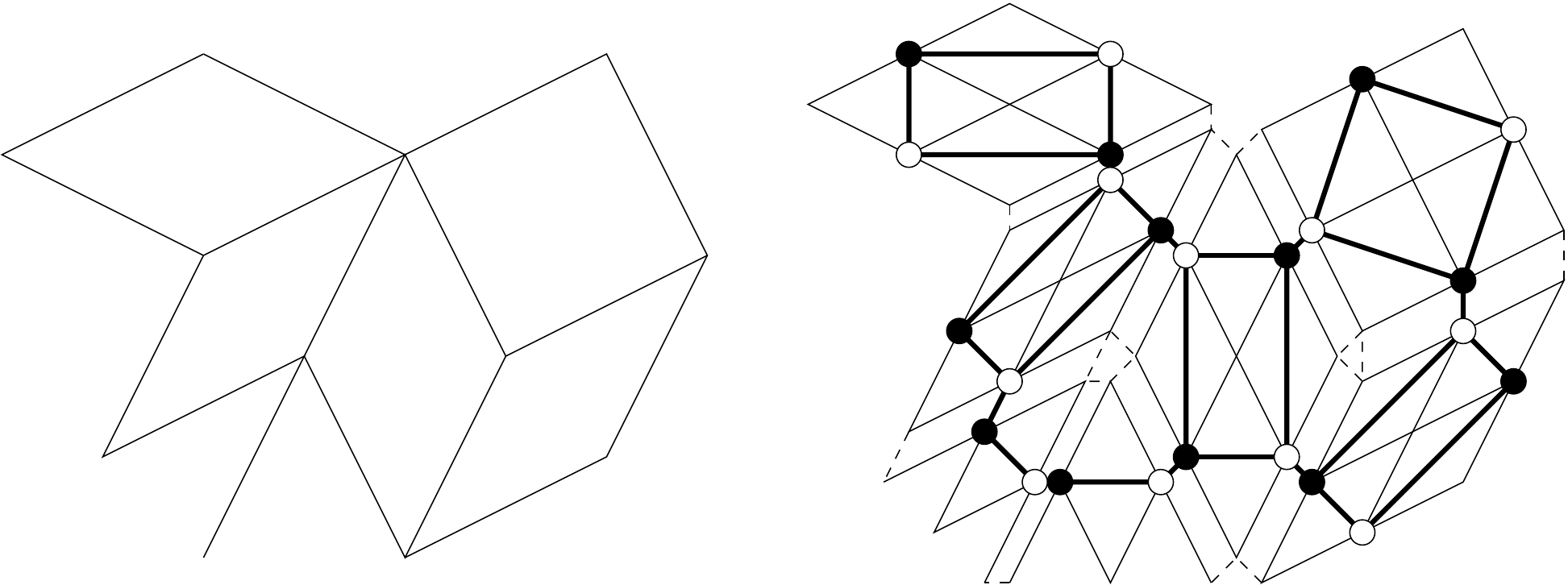}
\end{center}
\caption{a rhombus tiling $\dmd$; associated rhombus tiling $\dmd'$ and bipartite isoradial $C$ (with flat rhombi blown up for legibility; to be glued back along dashed lines)}
\label{fig:isorad2}
\end{figure}

Let us consider $C$, which is up to now embedded rather arbitrarily. Each pair of a black $b$ and white $w$ vertex connected by a road in $C$ corresponds to a pair $v\in\Gamma,f\in\Gamma^\dg$, ie to an edge in $\dmd$; we choose the coloring of $C$ s.t. $({\vec vf},{\vec bw})$ is direct (as in Figure \ref{fig:kastorient}). Let us embed both $b$ and $w$ at the midpoint of this edge of $\dmd$ (this is not a proper graph, but is nonetheless convenient here). Then $C$ is itself isoradial, and corresponds to the diamond graph $\dmd'$ obtained from $\dmd$ by dividing each face of $\dmd$ in four isometric rhombi, and inserting a flat rhombus on each edge of $\dmd$. This is somewhat degenerate but preserves most of the isoradial machinery (in particular the all-important discrete exponential functions, \cite{Mercat,Ken_isoradial,SmiCHe_isoradial}). See Figure \ref{fig:isorad2}.

Associated to any bipartite isoradial graph $C_0$ (with diamond graph $\dmd_0$) is a Dirac operator $K$ introduced in \cite{Ken_isoradial}. It is a complex weighted adjacency matrix with matrix elements: $K(w,b)=i(u-v)$ if $w,v,b,f$ are vertices on a face of $\dmd_0$ listed in counterclockwise order, with $b\in V^B_{C_0}$, $w\in V^W_{C_0}$ (adjacent black and white vertices). This defines $K:\C^{V^B_{C_0}}\rightarrow\C^{V^W_{C_0}}$. 

In the special case of $\dmd'$ and $C$, for each pair $v\in\Gamma,f\in\Gamma^\dg$, we have chosen an angle $\nu(vf)$ s.t. $f-v=|f-v|e^{2i\nu(vf)}$. The pair $(v,f)$ corresponds to a pair $b,w$ of vertices in $C$; set $\nu(b)=\nu(w)=\nu(vf)$. Then:
$$K(w,b)=ie^{i(\nu(w)+\nu(b))} K_C(w,b)$$
where $K_C$ is the (real) Kasteleyn operator we have been using and $K$ the critical Dirac operator. This is a simple direct check (distinguishing the cases of roads and streets), which informed the choice of Kasteleyn orientation on $C$.

\subsection{Fermionic variables}

For the reader's convenience, we include a discussion of fermionic variables and associated correlators in some of their various incarnations.

Following Kadanoff and Ceva \cite{KC_disorder}, we have considered a fermionic variable as the product of adjacent order and disorder variables: $\psi(vf)=\sigma(v)\mu(f)$ (which is invariant by Kramers-Wannier duality). This leaves an ambiguity on the sign of a correlator containing $\psi$, which is insufficient for some local computations. This ambiguity may be resolved in a couple of different ways; let us indicate one which is immediately compatible with the previous discussion, restricted for simplicity to the case of a pair fermionic correlator $\langle\psi(v_1f_1)\psi(v_2f_2)\rangle$ in a simply connected domain.

In order to assign a sign to this correlator, it is enough to specify a disorder line from $f_1$ to $f_2$, or rather an equivalence class of disorder lines. This may be done in terms of the choice of angles $\nu(v,f)$. In a slightly more general setting, one can consider an assignment $(v,f)\mapsto 2\nu(v,f)\in\R$ s.t. $\arg(f-v)=2\nu(v,f)\mod 2\pi$. Given $\nu(v_1,f_1)$ and a simple path $\gamma:v_1\rightarrow v_2$, one can define a natural parallel transport along $\gamma$ by the following local rules: there is a sequence $((v^k,f^k))_{k=0,\dots n}$ such that $(v^0,f^0)=(v_1,f_1)$ (resp. $(v^n,f^n)=(v_2,f_2)$) and for all $k$, $((v^k,f^k),(v^{k+1},f^{k+1}))$ is of one of two types: $((v,f),(v',f))$ with $v,v'$ consecutive vertices around $f$ in clockwise order; or $((v,f),(v,f'))$ with $f,f'$ consecutive points on $\gamma$ and $v$ the vertex on the lefthand side of $\vec{ff'}\in E_{\Gamma^\dg}$. We may define sequentially angles $\nu'(v^kf^k)$ so that $\nu'(v^0f^0)=\nu(v_0f_0)$ and the increment $\nu'(v^{k+1}f^{k+1})-\nu'(v^{k}f^k)$ is in $(-\pi,0)$ in the first case and in $(0,\pi)$ in the second case (with $\arg(f^k-v^k)=2\nu'(v^kf^k)\mod 2\pi$ for all $k$). Let us denote $T_{\gamma}\nu(v_1f_1)=\nu'(v^nf^n)$. It is easy to check that $T_{\gamma^r}T_\gamma\nu=\nu-\pi$ if $\gamma^r:f_2\rightarrow f_1$ is $\gamma$ taken with reverse orientation.

Given a choice of $\nu$'s, we may specify the sign of $\langle \psi(v_1f_1)\psi(v_2f_2)\rangle$ by evaluating this correlator with respect to a disorder line $\gamma:f_1\rightarrow f_2$ satisfying: $\nu(v_2f_2)=T_\gamma\nu(v_1f_1)\mod 2\pi$; in other words, if $\langle \psi(v_1f_1)\psi(v_2f_2)\rangle_\gamma$ is the evaluation of the correlator wrt $\gamma$, we may set
$$\langle\psi(v_1f_1)\psi(v_2f_2)\rangle=e^{i(\nu'(v_2f_2)-T_\gamma\nu(v_1f_1))}\langle \psi(v_1f_1)\psi(v_2f_2)\rangle_\gamma$$
Remark that with these conventions, the $\psi$'s are (as they should) anticommuting variables, viz. $\langle\psi(v_1f_1)\psi(v_2f_2)\rangle=-\langle\psi(v_2f_2)\psi(v_1f_1)\rangle$. 

Let us consider an Ising model on a simply-connected $\Gamma^\dg$ (spin variables are defined on faces) in its low-temperature expansion (a polygon on $\Gamma$). If $(v_if_i)$, $i=1,2$ are two vertex-face pairs, $\sigma^\dg(f_2)\sigma^\dg(f_1)=(-1)^n$ where $n$ is the number of polygon edges crossed by any path $\gamma:f_1\rightarrow f_2$ ($n\mod 2$ does not depend on the choice of $\gamma$, since the polygon has even degree at each $v\in\Gamma$). Disorder variables (which are here located on $V_\Gamma$) have the following low-temperature representation: a polygon with disorders at $v_1,v_2$ is a subgraph of $\Gamma$ with even degree at every $v\in V_\Gamma\setminus\{v_1,v_2\}$ and odd degree at $v_1,v_2$. Then the correlator $\langle\psi^\dg(f_1v_1)\psi^\dg(f_2v_2)\rangle$ may be represented in terms of polygons with disorders at $v_1,v_2$ and counted with a sign depending on the parity of the number of polygon edges crossed by $\gamma:f_1\rightarrow f_2$ (which does not cross a specified disorder line $v_1\rightarrow v_2$).

From the low-temperature expansion, one may move to the Fisher representation. From a polygon on $\Gamma$ with disorders at $v_1,\dots,v_{2n}$, one can define a matching on the Fisher graph $\Gamma_F$. Each polygon edge descends to a matched road on $\Gamma_F$; within each regular ($\neq v_1,\dots,v_{2n}$) city there are two choices to complete the matching. In each disordered city (corresponding to one of $v_1,\dots,v_{2n})$), remove an inner (A-type) vertex corresponding to the pair $(v_if_i)$. Then it is easy to check that the road matching may be extended in a unique way in the city with a monomer at $a_i\simeq(v_if_i)$. Consequently, in terms of the dimer configuration on $\Gamma_F$, the correlator $\langle\psi(v_1f_1)\psi(v_2f_2)\rangle$ is expressed (up to a multiplicative factor $2^2$ coming from the disordered decorations) as the partition function of matchings on $\Gamma_F\setminus\{v_1,\dots,v_{2n}\}$ counted with a sign depending on the parity of matched edges crossing a path $f_1\rightarrow f_2$. As is well-known, this is (up to sign) the inverting kernel of $K_F$ evaluated at the pair $a_1,a_2$: $K_F^{-1}(a_1,a_2)$.

As discussed earlier, $K_F^{-1}$ may be expressed in terms of $K_C^{-1}$ by linear algebra manipulations. Alternatively, by bosonization the correlator $\langle\psi(v_1f_1)\psi(v_2f_2)\rangle$ may be expressed in terms of magnetic/electric dimer correlators on $K_C$, which in turn are readily identified in terms of the inverting kernel $K_C^{-1}$. 

A variant of the (low temperature) polygon representation uses edge defects: every vertex has even degree, but in some prescribed edges one half-edge is present in the polygon and one half-edge is absent. Consider the partition function of these polygons with defects at $e_1=(v_1v_1'),e_2=(v_2v'_2)\in E_\Gamma$, where as before configuration are counted with a sign depending on parity; this is essentially the spin Ising observable used in \cite{SmiChe_Ising}.
 By erasing the defective half-edges, one obtains a polygon with one vertex defect at $v_1$ or $v'_1$ and one vertex defect at $v_2$ or $v'_2$, where the edges $e_1,e_2$ are vacant. It is then easy to see that this observable may be expressed as a linear combination of pair fermionic correlators.

As explained in \cite{NieKno_potts}, the correlator $\langle\psi(v_1f_1)\psi(v_2f_2)\rangle$ has a natural representation in terms of the random-cluster representation of the model (or rather the associated fully-packed loop representation), in the simply connected, wired boundary case. In the case of a simply connected domain with a wired and a free boundary arc, if $(v_1f_1)$ is taken ``in the bulk" and $(v_2f_2)$ is at one of the boundary condition change point on the boundary, one obtains the FK Ising observable of \cite{Smi_ICM} (up to normalization).

An important property of fermionic correlators $\langle\psi(vf)X\rangle$ is that they satisfy local linear relations. In \cite{Dotsenko_Ising}, it is shown (on the square lattice, see \cite{Mercat} on isoradial lattices) how to obtain a relation on the $\langle\psi(vf_i)X\rangle$'s, where the $f_i$'s are the faces adjacent to $v$ (and by duality on the $\langle\psi(v_jf)\dots\rangle$, where the $v_j$'s are the corner of $f$). Here $X$ represents arbitrary (fixed) order and disorder insertions away from $vf$. Linear relations of this type are instrumental in determining the scaling limit of FK and spin interfaces \cite{Smi_ICM}.
 
These relations may be derived in any of the various representations of  fermionic correlators; let us sketch one possible way for the Kadanoff-Ceva representation. Consider:
$$F(vf)=\langle\psi^\dg(vf)X\rangle$$
where $X$ represents order and disorder insertions away from $vf$ (say more than one lattice spacing away; order variables are on faces in order to simplify the comparison with the Fisher representation). Let $(vv')$ be an edge of $\Gamma$ and $(ff')$ be the dual (oriented) edge in $\Gamma^\dg$; we want to relate the four values $F(vf),F(v'f),F(vf'),F(v'f')$. Let us consider partial partition functions ${\mc Z}_{\pm,\pm}$ where 
$${\mc Z}_{\eps,\eps'}=\sum_{\stackrel{(\sigma_u)_{u\in V_{\Gamma^\dg}}}{\sigma(f)=\eps,\sigma(f')=\eps'}}X(\sigma)\prod_{\stackrel{e=(uu')\in E_{\Gamma^\dg}}{\sigma(u)\neq\sigma(u'),e\neq(ff')}}w_e$$
Remark that if we choose disorder lines, a pair of disorder variables can be represented as a (nonlocal) random variable. For definiteness, assume that $(vv')$ goes from left to right, $(ff')$ from bottom to top; we choose a disorder line $\gamma$ from $v$ to $v''$ (another disorder variable included in $X$) which does not contain $(vv')$; and a disorder line $\gamma'$ starting from $v'$ obtained by concatenating $(v'v)$ and $\gamma$. We use a local choice of $\nu$'s as in Figure \ref{fig:kastorient} ($F(vf)=\pm\langle\psi^\dg(vf)X\rangle_\gamma$, where $\pm$ depends on the choice of $\nu$'s and $\gamma$). Then
$$\begin{pmatrix} 
      -F(vf')\\
      F(vf)\\
      F(v'f)\\
      F(v'f')
   \end{pmatrix}
=   \pm\begin{pmatrix} 
      1&-1&-w&w\\
      1&-1&w&-w\\
      w&-w&1&-1\\
       w&-w&-1&1    
   \end{pmatrix}
\begin{pmatrix} 
      {\mc Z}_{++}\\
      {\mc Z}_{--}\\
      {\mc Z}_{+-}\\
      {\mc Z}_{-+}
   \end{pmatrix} 
$$
where $w=w_{(vv')}=\tan(\theta/2)$. Consequently:
$$
\begin{pmatrix} 
-1+w^2&1+w^2&-2w&0\\
-2w&0&1-w^2&1+w^2
   \end{pmatrix}
\begin{pmatrix} 
      -F(vf')\\
      F(vf)\\
      F(v'f)\\
      F(v'f')
   \end{pmatrix}=\begin{pmatrix} 
     0\\
     0
   \end{pmatrix}
$$
Each pair $(vf)$ corresponds to a black vertex of $C$; hence we may see $F$ as an element of $\R^{V_C^B}$, and we have just checked (see \eqref{eq:Kholom}) that $K_C^{BW}F=0$ (away from other insertions, which create monodromies and/or poles).

More generally, if $f\in\R^{V^B_C}$ and $K_C^{BW}f=0$ in some region, then around an edge $e\in E_\Gamma$ as in \eqref{eq:Kholom}, we have:
$$ 
\begin{pmatrix} 
     \cos(\theta)&1&-\sin(\theta)&0\\
     \sin(\theta)&0&\cos(\theta)&1
   \end{pmatrix}
\begin{pmatrix} 
     f(a_1)\\ f(a_2)\\ f(a_3)\\ f(a_4)   \end{pmatrix}
=\begin{pmatrix} 
     0\\ 0
   \end{pmatrix}
$$
ie $(f(a_1),\dots,f(a_4))$ lies in a two-dimensional real vector space,  which may be identified with $\C$ via:
$$z_e\mapsto \left(\Re(z_e e^{-i\nu(a_1)}),\dots,\Re(z_e e^{-i\nu(a_4)})\right)$$
Here $\theta\in (0,\pi/2)=\frac 12\arg\frac{f'-v}{f-v}$ and $(\nu(a_1),\dots,\nu(a_4))=(\frac\theta 2,\pi-\frac\theta 2,\frac\pi 2+\frac\theta 2,\frac{3\pi}2-\frac\theta 2)$. Hence associated to $f$ we have $g:E_\Gamma\rightarrow\C$ such that if $e,e'\in E_\Gamma$ are the two edges of $\Gamma$ corresponding to $a\in V^B_C$, we have the $\R$-linear relation:
$$\Re(z_{e'}e^{-i\nu(a)})=f(a)=\Re(z_ee^{-i\nu(a)})$$
which is (up to a rotation of the $z$'s) the notion of {\em $S$-holomorphicity} (\cite{SmiChe_Ising}).

\subsection{Quadratic relations}

In \cite{Perk_quad}, Perk identifies quadratic identities between spin correlators, which may be interpreted as discrete time Toda equations. Let us briefly indicate how to recover such identities (on a general graph) from bosonization. 

We consider dual Ising models on (simply-connected) $\Gamma,\Gamma^\dg$; for $i=1,2$, $e_i=(v_iv'_i)\in E_\Gamma$ and $e_i^\dg=(f_if'_i)$. Let us start from
$$\langle\sigma(v_1)\sigma(v_2)\rangle_\Gamma
\langle\sigma(v'_1)\sigma(v'_2)\rangle_\Gamma
-\langle\sigma(v_1)\sigma(v'_2)\rangle_\Gamma
\langle\sigma(v'_1)\sigma(v_2)\rangle_\Gamma$$
which by Kramers-Wannier duality is
$$c\left(\langle\sigma(v_1)\sigma(v_2)\rangle_\Gamma
\langle\mu^\dg(v'_1)\mu^\dg(v'_2)\rangle_{\Gamma^\dg}
-\langle\sigma(v_1)\sigma(v'_2)\rangle_\Gamma
\langle\mu^\dg(v'_1)\mu^\dg(v_2)\rangle_{\Gamma^\dg}
\right)
$$
where throughout $c$ is an (explicit) product of local factors (independent of insertions). By bosonization, this is rewritten as:
$$c\left(\langle\mu(v_1)\mu(v_2)\mu(v'_1)\mu(v'_2)\sigma(v'_1)\sigma(v'_2)\rangle_{6V}
-\langle\mu(v_1)\mu(v'_2)\mu(v'_1)\mu(v_2)\sigma(v'_1)\sigma(v_2)\rangle_{6V}
\right)$$
The edges $e_1,e_2\in E_\Gamma$ are identified with a pair of vertices on the 6V graph $M$. In both cases, we have a 6V configuration with disorder lines from $v_1$ to $v'_1$ and $v_2$ to $v'_2$, which may be drawn so as to intersect two consecutive edges of $M$ around $e_1,e_2$. A direct examination shows that terms corresponding to a sink and a source around $e_1$ and a sink and a source around $e_2$ are counted with the same sign in both correlators; this leaves terms corresponding to two sinks around $e_1$ and two sources around $e_2$ or vice-versa. Remark that having two sinks on edges around $e_1$ forces the 6V type of $e_1$; and that such a pair of edge defects may be alternatively represented by a vertex defect, where the arrow configuration around $v_1$ is an 8V sink.

Hence the above combination of spin correlations may be expressed as a 6V correlator with a pair of (magnetic) 8V defects at $e_1,e_2$ counted with a sign (electric) depending on the height variation from, say, $v'_1$ to $v_1$. This is plainly symmetric in $\Gamma\leftrightarrow\Gamma^\dg$ and consequently:
\begin{align*}
\langle\sigma(v_1)\sigma(v_2)\rangle_\Gamma
\langle\sigma(v'_1)\sigma(v'_2)\rangle_\Gamma
-\langle\sigma(v_1)\sigma(v'_2)\rangle_\Gamma
\langle\sigma(v'_1)\sigma(v_2)\rangle_\Gamma=\\
c\left(
\langle\sigma^\dg(f_1)\sigma^\dg(f_2)\rangle_{\Gamma^\dg}
\langle\sigma^\dg(f'_1)\sigma^\dg(f'_2)\rangle_{\Gamma^\dg}
-\langle\sigma^\dg(f_1)\sigma^\dg(f'_2)\rangle_{\Gamma^\dg}
\langle\sigma^\dg(f'_1)\sigma^\dg(f_2)\rangle_{\Gamma^\dg}\right)
\end{align*}
where $c$ is an explicit local factor depending on the choice of normalization.

\subsection{Consequences for dimers}

Bosonization allows to transfer questions about Ising correlations into problems on (bipartite) dimer correlator. The latter are typically analytically easier to handle and come with a natural free field interpretation. Let us point out however that some Ising results may be transferred back to dimers; we shall mention a few of these.

For the dimer height function on the square lattice, we know \cite{Ken_domino_GFF} that the scaling limit is the free field, which in the plane has the distributional invariance $\phi\leftrightarrow -\phi$. However at the discrete level this is not immediately apparent. An easy way to see it is to map the dimer configuration to a 6V configuration, use the arrow reversal of the latter, and map back to a dimer configuration. The mapping is (almost) deterministic and (almost) involutive, with ambiguities coming only from type 5-6 vertices (Figure \ref{fig:6Vdimer}).

In \cite{Dub_tors}, the asymptotics of electric correlators of type $\langle\exp(i\alpha(\phi(y)-\phi(x)))\rangle$ are evaluated for a charge $\alpha\in (-\frac 12,\frac 12)$, the limiting case $\alpha=\frac 12$ appearing as rather delicate. On the planar square lattice $\Z^2$ (with critical Ising weights), we may write
$$\E_{\Z^2}\left(\prod_{i=1}^{2n}\sigma(v_i)\right)\E_{\Z^2+u}\left(\prod_{i=1}^{2n}\sigma(v_i+u)\right)=
\E_{\Z^2}\left(\prod_{i=1}^{2n}\sigma(v_i)\right)^2
$$
where $u=\frac{1+i}2$ and $(\Z^2)^\dg\simeq u+\Z^2$. Using bosonization on both sides yields:
$$\E_{\rm dimer}\left(\prod_{i=1}^{2n}({\mc O}_1+{\mc O}_{-1})(v_i+\frac u2)\right)=
\E_{\rm dimer}\left(\prod_{i=1}^{2n}\cos(\phi(v_i)/2)\right)
$$
where the dimer graph is $\frac 12(u+\Z^2)$. This is an exact identity between electric and magnetic dimer correlators, and is a discrete version of (a by-product of) $T$-duality for the limiting free field. Since the monomer correlators are worked out \cite{Dub_tors}, this also allows (for the square lattice) to evaluate electric correlators with half-integer electric charge.

Lastly, let us consider the case of a subgraph $C$ of $\delta\Z^2$ (carrying dimers) corresponding to an Ising model on $\delta(\frac u2+2\Z^2)$ with wired (or dually free) boundary conditions along (a $\delta$-approximation of) a fixed smooth closed loop $\gamma$. On the one hand, the inverse Kasteleyn matrix $K_C^{-1}(.,.)$ may be identified as a (properly normalized) pair fermion correlator. On the other hand, such pair correlations (under a different representation, see above) are analyzed in \cite{HonSmi_energ}, building in particular on the discrete integration of squares of $S$-holomorphic functions. This shows that the asymptotics of $K_C^{-1}$ derived above (in an a rather pedestrian fashion) for polygonal domains actually holds for a much larger class of (simply-connected) domains. Let us repeat that this choice of boundary corresponds to a flat boundary height for dimers. Then convergence of the height field to a free field with constant (Dirichlet) boundary conditions follows easily from Sects 3-5 in \cite{Dub_tors}.

{\bf Acknowledgments.} I wish to thank Cl\'ement Hongler for very interesting conversations during the preparation of this article.

\bibliographystyle{abbrv}
\bibliography{biblio}

-----------------------

\noindent Columbia University\\
Department of Mathematics\\
2990 Broadway\\
New York NY 10027

\end{document}